\documentclass[sn-mathphys]{sn-jnl}% Math and Physical Sciences Reference Style
\usepackage{amssymb,amsmath,mathrsfs, amsthm}

\usepackage{cases}

\jyear{2021}%

\theoremstyle{thmstyleone}%
\numberwithin{equation}{section}
\newtheorem{theorem}{Theorem}[section]%  meant for continuous numbers
\newtheorem{question}{Question}[section]%

\newtheorem{proposition}[theorem]{Proposition}%
\newtheorem{corollary}{Corollary}[section]
\theoremstyle{thmstyletwo}%
\newtheorem{remark}[theorem]{Remark}%
\newtheorem{lemma}{Lemma}[section]%

\theoremstyle{thmstylethree}%
\newtheorem{definition}[theorem]{Definition}%

\begin{document}

\title[On strongly pseudoconvex K\"ahler--Berwald manifolds]{On strongly pseudoconvex K\"ahler--Berwald manifolds}

\author[]{\fnm{Xiaoshu} \sur{Ge}}\email{gxs@amss.ac.cn}

\author[]{\fnm{Xiangyu} \sur{Zhou}}\email{xyzhou@amss.ac.cn}


\affil{\orgdiv{Institute of Mathematics, Academy of Mathematics and Systems Science}, \orgname{Chinese Academy of Sciences},  \city{Beijing }, \postcode{100190}, \country{China}}

%\affil[2]{\orgdiv{Institute of Mathematics, Academy of Mathematics and Systems Science}, \orgname{Chinese Academy of Sciences}, \orgaddress{ \city{Beijing}, \postcode{100190}, \country{China}}}


\abstract{
In this paper, we establish the equivalence between the K\"ahler--Finsler and K\"ahler--Berwald conditions for complex Finsler metrics. That is, a K\"ahler--Finsler metric always satisfies the complex Berwald condition: its Chern--Finsler nonlinear connection reduces to the Chern connection of a K\"ahler metric. We also provide, for the first time, smooth weakly K\"ahler--Finsler examples which are not K\"ahler--Berwald; and propose a condition that is formally stronger than weakly K\"ahler--Finsler, with equivalence holding in complex dimension two.
Then we prove the de Rham decomposition and uniformization theorems for K\"ahler--Berwald manifolds, and show that every strongly pseudoconvex locally symmetric complex Finsler metric is both K\"ahler--Berwald and complex Einstein--Finsler.
Additionally, we establish several rigidity theorems for the Kobayashi and Carath\'eodory metrics: $F_K$ is K\"ahler--Berwald if and only if it is K\"ahler;
if the manifold is simply connected, then $F_C$ is K\"ahler--Berwald if and only if it is K\"ahler.
}

%%================================%%
%% Sample for structured abstract %%
%%================================%%




\maketitle

\section{Introduction and main results}\label{int-res}

In this paper, we obtain the coincidence between the K\"ahler--Finsler and K\"ahler--Berwald conditions, and
propose a formally stronger condition than weakly K\"ahler--Finsler. We prove that the new condition coincides with weakly K\"ahler--Finsler in complex dimension two, and their coincidence has also been established in \cite{Ge-wKF}.
Then we establish the de Rham decomposition, uniformization and rigidity theorems in the K\"ahler--Berwald setting.
We also prove that the $h\bar h$-part curvature of each locally irreducible de Rham factor of a K\"ahler--Berwald manifold $(M,F)$ is just the restriction of the $h\bar h$-part curvature of $F$.
Thus many results in K\"ahler geometry can be generalized to the K\"ahler--Berwald setting, for example, the generalized Frankel conjecture, which is verified in this paper.

\subsection{Generalizations of the K\"ahler condition in complex Finsler geometry}\label{subsec1.0}

 Let $M$ be a complex manifold of complex dimension $n$. A complex Finsler metric on $M$ is defined as follows.
\begin{definition}[\cite{AP}]
A smooth complex Finsler metric on $M$ is an upper
semi-continuous function $F:T^{1,0}M\rightarrow [0,+\infty)$
such that
\begin{itemize}
  \item [1)] $G:=F^2$ is smooth on 	
 $\tilde{M}$, where $\tilde{M}:=T^{1,0}M\setminus\{\text{zero section}\}$;
  \item [2)] $G>0$ on $\tilde{M}$;
  \item [3)]  $F(z,\lambda v)=\vert\lambda\vert F(z,v)$ for all $(z,v)\in \tilde M$ and $\lambda\in\mathbb{C}$.
\end{itemize}
And $F$ is called
\begin{itemize}
  \item [4)] strongly pseudoconvex if furthermore the Levi matrix $(G_{i\overline{j}}):=\left(\frac{\partial^2G}{\partial
v^i\partial\overline{v^j}}\right)$ is positive definite on $\tilde{M}$, which is equivalent to requiring that all the $F$-indicatrices $ I_F(p):=\{v\in T^{1,0}_pM:F(p,v)<1\}$ are strongly pseudoconvex domains;
\item [5)] strongly convex if in addition the real Hessian matrix $\left(\frac{\partial^2G}{\partial
u^i\partial{u^j}}\right)$ is positive definite on $\tilde{M}$ when we identify
$z^j = x^j + \sqrt{-1} x^{n+j}$ and $v^j = u^j + \sqrt{-1} u^{n+j}$, and take $G(x;u):= G(z,v)$. This is equivalent to the condition that all $F$-indicatrices are strongly convex domains.
\end{itemize}
\end{definition}
\begin{remark}
For a strongly pseudoconvex complex Finsler metric $F$, $G$ is smooth on $T^{1,0}M$ if and only if $G$ is a smooth Hermitian (quadratic) metric \cite{AP}.
\end{remark}

%According to S.-S. Chern \cite{Chern}, Finsler geometry is just Riemannian geometry without the quadratic restriction.
For a strongly pseudoconvex complex Finsler metric $F$, applying the Euler formula for homogeneous maps to
\[G(z,v):=F^2(z;\lambda v)=\vert\lambda\vert^2G(z,v),\quad\forall \, (z,v)\in \tilde{M},\,\lambda\in \mathbb C,\]
 one has
\begin{equation}
G(z,v)=G_{i\bar j}(z,v)v^i\bar v^j,\quad\forall \,(z,v)\in \tilde{M},\label{G-fund}
\end{equation}
where $G_{i\bar j}:=\frac{\partial^2G}{\partial v^i\partial \bar v^j}.$
By \eqref{G-fund}, $G$ is induced from a Hermitian metric if and only if
\[G_{i\bar j}(z,v)=G_{i\bar j}(z),\quad \forall \,(z,v)\in \tilde M.\]
This leads to an essential difference between non-Hermitian (quadratic) complex Finsler metrics and Hermitian (quadratic) ones in their connection coefficients, hence also in torsion and curvature: the former depend on tangent vectors, while the latter do not.
In the adapted local coframe $\{dz^k,\psi^k\}_{k=1}^n$, dual of the local frame $\{\partial_k,\delta_k\}_{k=1}^n$ for the tangent bundle, as provided in Section \ref{Sec-cpx-F}, the torsion forms of the Chern--Finsler connection of $F$ contain two parts \cite{AP93}:
\begin{align}
 \theta^i
 &=-dz^k\wedge\omega_k^i=\frac12(\Gamma^i_{j;k}-\Gamma^i_{k;j})dz^k\wedge dz^j+\Gamma^i_{jk}\psi^k\wedge dz^j; \label{tor20}\\
 \dot\theta^i&=\partial \psi^i-\psi^k\wedge\omega_k^i=0;\\
 \tau^i&=\bar \partial\psi^i=-\delta_{\bar j}\Gamma^i_k\,dz^k\wedge d\bar z^j
          -\Gamma^i_{\bar j;k}\,dz^k\wedge\bar\psi^j. \label{tor11}
\end{align}
Using the vanishing of the horizontal $(2,0)$-part torsion, $\Gamma^a_{j;k}-\Gamma^a_{k;j}=0$, H. Rund introduced the
first Finsler-type generalization of the K\"ahler condition in 1967 \cite{Rund}, which is renamed as strongly K\"ahler--Finsler after many years by  Abate and Patrizio in 1994 \cite{AP93}.
Also in \cite{AP93}, Abate--Patrizio obtained the following rigidity result.
\begin{proposition}[{\cite{AP93}}]
 Let $F$ be a strongly pseudoconvex complex Finsler metric on a complex manifold $M$. Then:
\begin{itemize}
  \item [1)] the mixed part of the (2,0)-torsion $\Gamma^i_{jk}$ vanishes if and only if $F$ comes from a Hermitian metric;
  \item [2)] $\theta^i$ vanishes if and only if $F$ comes from a K\"ahler metric.
\end{itemize}
\end{proposition}
Then, Abate--Patrizio introduced the K\"ahler--Finsler and weakly K\"ahler--Finsler notions \cite{AP93}. In that work, they employed the weakly K\"ahler--Finsler condition in the first variation of the length integral to proved the local existence and uniqueness of geodesics for strongly convex metrics. They also established the second variation formula of the length integral under the K\"ahler--Finsler condition, and subsequently referred to the condition introduced by H. Rund as strongly K\"ahler--Finsler.
The complex Berwald condition was introduced by Aikou in 1991, who also investigated the relationship between complex Berwald metrics and their induced Hermitian ones \cite{Ai91}.
\begin{definition}\label{KF-KB-wKF}
A strongly pseudoconvex complex Finsler metric $F$ is called
\begin{itemize}
\item strongly K\"ahler--Finsler (SKF) if
\begin{equation}
\Gamma_{j;i}^{l}(z,v) = \Gamma_{i;j}^{l}(z,v), \quad
\forall \,(z,v) \in \tilde{M};
\label{Kahler}
\end{equation}
\item K\"ahler--Finsler (KF) if
\begin{equation}
(\Gamma_{j;i}^{l}(z,v) - \Gamma_{i;j}^{l}(z,v))v^j=0, \quad
\forall \,(z,v) \in \tilde{M};
\label{KahlerFinsler}
\end{equation}
\item weakly K\"ahler--Finsler (WKF) if
\begin{equation}
G_l(\Gamma_{j;i}^{l}(z,v)-\Gamma_{i;j}^{l}(z,v))v^j=0, \quad
\forall \,(z,v) \in \tilde{M};
\label{w-Kahler}
\end{equation}
\item complex Berwald (B) if
\begin{equation}
\Gamma_{i;j}^{l}(z,v) = \Gamma_{i;j}^{l}(z), \quad
\forall \,(z,v) \in \tilde{M};
\label{Berwald}
\end{equation}
\item K\"ahler--Berwald (KB) if both \eqref{KahlerFinsler} and \eqref{Berwald} hold.
\end{itemize}
\end{definition}
There is another condition, formally stronger than weakly K\"ahler--Finsler (SWKF), we do not find its appearance elsewhere, defined as follows
\begin{equation}
G_l(\Gamma_{j;i}^{l}(z,v)-\Gamma_{i;j}^{l}(z,v))=0, \quad
\forall \, (z,v) \in \tilde{M}.
\label{new-w-Kahler}
\end{equation}
By definition, the following relationships hold:
\[SKF\Rightarrow  KF\Rightarrow WKF, \, SKF\Rightarrow  SWKF\Rightarrow WKF,  \, KB\Leftrightarrow KF+B.\]
In 2009, B. Chen and Y. Shen demonstrated the following result \cite{CS}:
\[SKF\Leftrightarrow KF.\]
We obtain here, in case of complex dimension two (cf. Proposition \ref{new-KF-th}),
\[SWKF\Leftrightarrow WKF.\]
The first-named author proved the coincidence between  $SWKF$ and $WKF$ for general dimensions later, then he provided some smooth weakly K\"ahler--Finsler examples which are not K\"ahler--Finsler \cite{Ge-wKF,Ge-KF}.
Our rigidity theorem (cf. Theorem \ref{0Kob-KB}) implies that a weakly K\"ahler--Finsler type Kobayashi metric is never K\"ahler--Finsler unless it is K\"ahler. Hence the following strict relationships are clear:
\[KF\Rightarrow WKF,\quad KF\nLeftarrow WKF.\]
In \cite{XZ}, H. Xia and C. Zhong proved
\begin{equation}
KB\Leftrightarrow WKF+B.\label{XZ-B-wK-K}
\end{equation}
Then, H. Xia and C. Zhong asked whether there exists K\"ahler--Finsler metric which is not K\"ahler--Berwald \cite{XZ}.
We obtain a negative answer to their question by the following theorem.
\begin{theorem}[cf. Theorem \ref{thm:KF=KB}]\label{thm:KF=KB1}
 K\"ahler--Finsler metrics are K\"ahler--Berwald.
\end{theorem}
In summary, the following relationships hold:
\begin{equation}
KB\Leftrightarrow SKF\Leftrightarrow KF\Leftrightarrow WKF
+B\Rightarrow WKF\Leftrightarrow SWKF,\quad SWKF \nRightarrow KF.\label{relation}
\end{equation}

\begin{remark}
The following results also imply important relationships between these Finsler version generalizations and the Hermitian K\"ahler condition:
\begin{itemize}
  \item when $G$ is induced from a Hermitian metric, $\Gamma_{i;j}^{l}(z,v) = \Gamma_{i;j}^{l}(z)$ just gives the Chern connection. In this circumstance, both $KB$ and $WKF$ reduce to the K\"ahler condition;
  \item  in 2011,  Aikou proved that the base manifold of every K\"ahler--Finsler manifold is K\"ahler \cite{Ai11}. Recently, the first-named author proved that the base manifold of every weakly K\"ahler--Finsler manifold is also K\"ahler \cite{Ge-KF}.
\end{itemize}
\end{remark}
For the sake of convenience, most of the rest results in this paper are presented in the K\"ahler--Berwald (which is currently identical to K\"ahler--Finsler) setting; and we shall use the K\"ahler--Finsler and strongly K\"ahhler--Finsler conditions interchangeably in the following.

\subsection{Locally symmetric complex Finsler spaces}\label{subsec1.2}
Riemannian (resp. Hermitian) symmetric spaces are classical generalizations of space forms, which are close connected with Berwald (resp. K\"ahler--Berwald) spaces. 	
%Symmetric spaces were initiated by \'E. Cartan in 1926 \cite{Cart}, and classified by him in terms of the classical and exceptional simple Lie groups \cite{Cart36}.
%Chinese mathematicians L.-G. Hua and Q.-K. Lu also made great contributions to geometric functional theories of symmetric bounded domains and symmetric manifolds \cite{Hua,Lu56,Lu57,Lu58,Lu}.
%For general theories about symmetric spaces, one can refer to Helgason's book \cite{Heg}. For homogeneous and symmetric real Finsler spaces, one can refer to the book \cite{Deng2} written by S. Q. Deng.
By Berger and Simons' results about torsion-free linear connection \cite{Ber55,Sim62}, the holonomy group of an irreducible manifold $M$ acts transitively on the unit sphere of each tangent space unless $M$ is a symmetric space of rank $\geq2$.
On the other hand, Aikou proved that every K\"ahler--Berwald metric $F$ is associated with a K\"ahler metric $h$, such that the Chern connection of $h$ coincides with the Chern--Finsler nonlinear connection of $F$ in \cite{Ai95} (cf. Theorem \ref{Ai-thm}).
In order that the Chern--Finsler nonlinear connection of a non-Hermitian K\"ahler--Berwlad metric coincides with the Chern connection of a associated K\"ahler metric, it's natural to study on symmetric spaces.
Indeed, in \cite{Deng}, S. Deng and Z. Hou classified all irreducible symmetric complex Finsler spaces, which proved to be symmetric Hermitian spaces as underlying manifolds.
%Analogous to the Hermitian cases, symmetric complex Finsler manifolds are locally symmetric, and locally symmetric complex Finsler manifolds are homogeneous.
%Even though symmetric complex Finsler manifolds are Hermitian symmetric manifolds as underlying spaces \cite{Deng}, they may have different geometries.
They also proved that complete locally symmetric real Finsler spaces are real Berwald \cite{Deng07}.
 Thus there is a natural question:
\begin{question}\label{ques-symm-KB}
Whether strongly pseudoconvex (locally) symmetric complex Finsler metrics are K\"ahler--Berwald?
\end{question}
The answer to Question \ref{ques-symm-KB} is positive in many known special cases \cite{Zh2,Lin1,GZ,CGZ,Zh3},
%In \cite{Zh2}, C. Zhong proved that on $\mathbb B_n$, any $\operatorname{Aut}(\mathbb B_n)$-invariant (hence symmetric) strongly pseudoconvex complex Finsler metric must be a constant multiple of the Poincar\'e-Bergman metric, while on the polydisc $P_n \subset \mathbb{C}^n$ ($n\geq 2$), there exist infinitely many complete $\operatorname{Aut}(P_n)$-invariant strongly pseudoconvex symmetric K\"ahler--Berwald metrics which are non-Hermitian quadratic. Subsequently, S. Lin and C. Zhong \cite{Lin1} proved that all holomorphically invariant, strongly pseudoconvex complex Finsler metrics on the unit polydisc are K\"ahler--Berwald metrics. In \cite{GZ}, the first-named author and C. Zhong constructed infinitely many holomorphically invariant metrics on each of the four classical domains which were proved to be K\"ahler--Berwald. %In \cite{CGZ}, P. Cao, the first author and C. Zhong proved that on the complex Grassmann manifold $U(m+n)/U(m)\times U(n)$, every $U(m+n)$-invariant strongly pseudoconvex complex Finsler metric is a K\"ahler--Berwald metric.
%Furthermore, C. Zhong \cite{Zh3} proved that all holomorphically invariant strongly pseudoconvex complex Finsler metrics on the four classes of the classical domains are K\"ahler--Berwald.
where the proofs rely deeply on the explicit expressions of the holomorphic automorphism groups and the invariance condition of the metric under biholomorphisms.
%In order to obtain the Chern--Finsler connection, it is needed to calculate the inverse metrix of the metric metrix, which is difficult in general and overcome by special isometrical biholomorphic mappings there.
%However, for exceptional classical domains, more general reducible symmetric complex manifolds, and locally symmetric complex manifolds, it is difficult to express all holomorphic invariant complex Finsler metrics explicitly.
%What's more, a symmetric complex Finsler metric maybe not invariant under the whole holomorphic automorphism group $Aut(M)$ of the manifold $M$, it may be invariant just invariant under a subgroup of  $Aut(M)$. To address the question in full generality, a new approach based on their common geometric properties is needed.
The fundamental property for a locally symmetric space $M$ is that: for every point $p\in M$, there exists a local neighborhood $U\subset M$ of $p$, such that the geodesic symmetry map, with $p$ as an isolated fixed point, is  holomorphic isometric.
Such holomorphically isometric geodesic symmetries with isolated fixed points are essential in answering Question \ref{ques-symm-KB}, which can be used to obtain the horizontal Chern--Finsler connection, thus to obtain the K\"ahler--Berwald conclusion. We obtain the positively answer to Question \ref{ques-symm-KB} in two ways (cf. Theorem \ref{KF=symm} and Theorem \ref{thm1}).
And we obtain the following result.
\begin{theorem}[cf.Theorem \ref{thm1}]\label{thm}
Every strongly pseudoconvex locally symmetric complex Finsler metric  $F$ is a K\"ahler--Berwald, complex Finsler--Einstein metric. And the corresponding Hermitian metrics whose Chern connection coincide with the Chern--Finsler nonlinear connection of $F$ are locally symmetric.
\end{theorem}

Using Theorem \ref{thm}, we further provide an equivalent characterization for locally symmetric complex Finsler metrics,
which is a generalization of the locally symmetric Hermitian version results.

\begin{theorem}[cf. Theorem \ref{pal-eq}]\label{1pal-eq}
A complete strongly pseudoconvex K\"ahler--Berwald metric $F$ is locally symmetric if only if
\[\nabla_H R^{\mathcal H} =0,\quad\forall \,H\in \mathcal H,\]
where $R^{\mathcal H}:=R_{j \bar s;k\bar l}dz^k\wedge d\overline{z^l}\otimes dz^j \otimes d\overline{z^s},$ and
$R_{j \bar s;k\bar l}$ is the $h\bar h$-part Chern--Finsler curvature defined by \eqref{HH-curv'}.
\end{theorem}

With similar methods used for the real symmetric Finsler case in \cite{Deng-Hou05}, and the result that strongly pseudoconvex locally symmetric complex Finsler metrics are K\"ahler--Berwald, we obtain the curvature formulas for strongly pseudoconvex locally symmetric complex Finsler metrics (cf. Proposition \ref{3curvature}).

\subsection{de Rham decomposition theorem}\label{subsec1.3}

The well known de Rham decomposition theorem due to G. de Rham %(cf. also Theorem 6.2 in Chapter IV of volume I \cite{Kob63})
says that \cite{DeRham}: a simply connected, complete Riemannian manifold $M$ is isometric to the direct product $M_0\times M_1\times\cdots\times M_n$, where $M_0$ is a Euclidean space (possibly of dimension 0) and $M_1,\cdots,M_n$ are all simply connected, complete, irreducible Riemannian manifolds. Such a decomposition is unique up to reordering. For K\"ahler manifolds, Hano and Matsushima obtained the following theorem \cite{Hano}.%, cf. also Theorem 8.1 in Chapter IX of volume II \cite{Kob63}).
\begin{theorem}[\cite{Hano}]\label{DEKahler}
A simply connected, complete K\"ahler manifold M is holomorphically isometric to the direct product
$M_0\times M_1\times\cdots\times M_n$, where $M_0$ is a complex Euclidean space (possibly of dimension 0), and $M_1,\cdots,M_n$ are all simply connected, complete, irreducible K\"ahler manifolds.
\end{theorem}

In \cite{szabo}, Szab\'o obtained the characterization for real Berwald manifolds, and generalized the de Rham decomposition theorem to the real Berwald setting.
%\begin{theorem}[\cite{szabo}]\label{reaDR}
%A connected (strongly convex real) Berwald space $(M,F,\nabla)$ must be one of the following four types:
%\begin{enumerate}
%\item [1)]  Riemann space;
%\item  [2)] locally Minkowski space;
%\item locally irreducible, locally symmetric, non-Riemannian Berwald space, of rank $r\geq 2$;
% \item  [3)] locally reducible and can be locally decomposed as a
%Descartes product of  Riemann space, locally Minkowski space, and locally irreducible, locally symmetric, non-Riemann Berwald spaces of rank $r\geq 2.$
% \end{enumerate}
%And if $(M,F,\nabla)$ is in addition simply connected and complete, it can be decomposed to the Descartes product of a Minkowskian space, simply connected complete irreducible Riemann space, and simply connected complete irreducible globally symmetric non-Riemannian Berwald space, of rank $r\geq 2$. Such a decomposition is unique up to an order.
%\end{theorem}
In \cite{Zh1}, C. Zhong proved that a strongly convex weakly K\"ahler--Finsler metric is a weakly complex Berwald metric if and only if it is a real Berwald metric. In \cite{Zh2},
C. Zhong obtained the following characterization and de Rham decomposition theorem for strongly convex K\"ahler--Berwald manifolds. A natural question is:
\begin{question}\label{ques-deRham}
 Does the de Rham decomposition theorem holds for strongly pseudoconvex K\"ahler--Berwald manifolds?
\end{question}
To answer Question \ref{ques-deRham}, we derive the relationship between affine translations and isometries for K\"ahler--Berwald metrics in Lemma \ref{aff-iso}, using similar methods applied to real Berwald manifolds in \cite{Deng2}.
%\begin{lemma}[cf. Lemma \ref{aff-iso}]\label{1all-iso}
% Let $(M_i,F_i)$, $i= 1,2,$ be two complex Berwald spaces, and $f$ an affine locally biholomorphism from $(M_1,F_1)$ onto $(M_2,F_2)$ with respect to the Chern--Finsler nonlinear connections of $F_1$ and $F_2$. If there exists $p\in M_1$ such that $df$ is a linear isometry from $T_p^{1,0}(M_1)$ onto $T_{\psi(p)}^{1,0}(M_2)$, then $f$ is a locally isometry.
%\end{lemma}
In \cite{Ai95}, Aikou associated to every strongly pseudoconvex K\"ahler--Berwald metric $F$ an induced K\"ahler metric $h$, such that the Chern connection of $h$ coincides with the Chern--Finsler nonlinear connection of $F$ (cf. Theorem \ref{Ai-thm}).
By Theorem \ref{DEKahler}, results about holonomy group by M. Berger and J. Simons \cite{Ber55,Sim62} (cf. Theorem \ref{Berger}),
and Lemma \ref{aff-iso}, we give an affirmative answer to Question \ref{ques-deRham}.
\begin{theorem}[cf. Theorem \ref{characKB}]\label{1characKB}
A connected, strongly pseudoconvex K\"ahler--Berwald manifold $(M,F,\nabla)$ must be one of the  following four types:
\begin{enumerate}
\item  [1)] K\"ahler space;
\item  [2)] locally complex Minkowski space;
\item  [3)] locally irreducible, locally symmetric, non-Hermitian, K\"ahler--Berwald manifold of rank $r \geq2$;
\item  [4)] locally reducible, and in this case $(M, F,\nabla)$ can be locally decomposed into a Descartes product of K\"ahler spaces, locally complex Minkowski spaces and locally irreducible, locally symmetric, non-Hermitian K\"ahler--Berwald manifolds of rank $r \geq2.$
    \end{enumerate}
\end{theorem}

\begin{remark}
Here the Descartes product for complex Berwald space is a generalization of general direct product of Hermitian manifolds.
The Berwald metric is not necessary a product metric as the Hermitian circumstance.
\end{remark}

\begin{theorem}[cf. Theorem \ref{KBDR}]\label{1KBDR}
A simply connected, complete, strongly pseudoconvex K\"ahler--Berwald manifold $(M,F,\nabla)$ can be decomposed into the Descartes product of a complex Minkowski space,  simply connected complete irreducible K\"ahler spaces and  simply connected non-Hermitian complete irreducible symmetric K\"ahler--Berwald manifolds of rank $r \geq2.$
\end{theorem}

The generalized Frankel conjecture, resolved by Mok \cite{Mok88}, also holds in the K\"ahler--Berwald setting.

\begin{theorem}[cf. Theorem \ref{KBFrankel4}]\label{KBFrankel}
Let $(M, F)$ be a compact, strongly pseudoconvex, K\"ahler--Berwald manifold of nonnegative holomorphic bisectional curvature.
Then the universal cover $(\widetilde{M}, \widetilde{F})$ of $(M, F)$ can be decomposed into the Descartes product with factors holomorphically isometric to:
\begin{enumerate}
  \item  [1)] complex Minkowski spaces;
  \item  [2)] complex projective spaces endowed with K\"ahler metrics of nonnegative holomorphic bisectional curvature;
  \item  [3)] irreducible compact-type symmetric complex Finsler spaces of rank $\geq2$.%, with each $F_l$ is unique on $M_l$ in the sense of Chern--Finsler nonlinear class.
\end{enumerate}
\end{theorem}
\begin{remark}
When $F$ is a K\"ahler metric, the decomposition of the Descarts product of $(M,F)$ reduces to the usual product of manifolds with usual product metric, and Theorem \ref{KBFrankel} reduces to the generalized Frankel conjecture for K\"ahler manifolds.
\end{remark}

\subsection{Uniformization and rigidity theorems}\label{subsec1.1}

The uniformization theorem of K\"ahler manifolds says that: a simply connected complete K\"ahler manifold of constant holomorphic sectional curvature must be holomorphically isometric (up to a multiplicative constant) to one of the complex Euclidean space $\mathbb C^n$ endowed with the complex Euclidean metric, the complex unit ball $\mathbb B^n$ endowed with the Bergman metric, and the complex projective space $\mathbb {P}^n$ endowed with the Fubini--Study metric.

Followed by a famous result of L. Lempert, on every smooth bounded strongly convex domain, the Kobayashi metric and the Carath\'eodory metric coincide, and they are complete, strongly convex, weakly K\"ahler--Finsler metrics of constant holomorphic sectional curvature $-4$ \cite{Lem81,Pat,AP3}.
Thus the uniformization theorem for the weakly K\"ahler--Finsler setting does not hold.

In \cite{AP3}, M. Abate and G. Patrizio obtained the uniformization theorem for K\"ahler--Finsler manifolds under the following curvature symmetry condition:
\begin{equation}
\langle \Omega(H,\bar\chi)\chi,\chi \rangle=\langle \Omega(\chi,\bar\chi)H,\chi \rangle,\forall \,H\in \mathcal H,\label{curv-sy}
\end{equation}
where $\chi:=v^i\delta_i$ is the radial horizonal vector field.

In \cite{CX-Yan}, X. Chen and R. Yan used a more symmetric condition for the curvature, called as the strictly condition in their paper:
\begin{equation}
\langle\bar\partial_{\mathcal H}\theta(H,\chi,\bar K),\chi\rangle=0,\quad \forall \,H,K\in \Gamma(\mathcal H).\label{CYKB}
\end{equation}
By \cite[Lemma 3.16]{AP2}, \eqref{CYKB} equivalents to:
\begin{equation}
\langle \Omega(H,\bar K)\chi,\chi \rangle=\langle \Omega(\chi,\bar K)H,\chi \rangle,\quad \forall \,H,K\in \mathcal H.\label{AP-curv-sy}
\end{equation}
Then they proved the uniformization theorem for K\"ahler--Berwald manifolds under the strictly condition \eqref{AP-curv-sy}.

In \cite{Ge2022}, the first-named author obtained for real Berwald strongly convex complex Finsler manifold $(M,F,J,\hat\nabla)$: the complex structure $J$ is parallel with respect to its Cartan nonlinear connecction if and only if $F$ is K\"ahler--Berwald.
Recently, W. Xiao and C. Zhong obtained for strongly convex complex Finsler metrics the coincidence between the horizontal parallelism of the complex structure with respect to its Cartan connection and the K\"ahler--Berwald condition. Based on this equivalent result, they obtained the uniformization theorem for K\"ahler--Berwald manifolds under the strongly convex assumption \cite{XZ26}.

In this paper, we establish the uniformization theorem for K\"ahler--Berwald menifolds without the strongly convex condition, and our proof is different from the above Finsler version results in \cite{Ai95,AP2,XZ26}.

\begin{theorem}[cf. Theorem \ref{uniformization}]\label{1uniformization}
Every simply connected, complete, K\"ahler--Berwald manifold of constant holomorphic sectional curvature $c$ is holomorphically isometric to
\begin{enumerate}
  \item [1)] $\mathbb B^n$ endowed with the Bergman metric if  $c=-4$;
  \item [2)] $\mathbb C^n$ endowed with a Minkowski metric  if $c=0$;
  \item [3)] $\mathbb {P}^n$ endowed with the Fubini--Study metric if $c=4$.
\end{enumerate}
\end{theorem}

%\begin{remark}
%The methods used here is different from the previous results in \cite{Ai95,AP2,XZ26} as presented at Section \ref{sub-cHSC}. Which contain the following main steps:
%\begin{enumerate}
%  \item we will use the de Rham decomposition theorem for strongly pseudoconvex K\"ahler--Berwald metrics established in this paper as Theorem \ref{1KBDR};
%\item we prove that a manifold under the setting of Theorem \ref{1uniformization} is either a complex Minkowski space or an irreducible manifold, and in the later case the metric must be K\"ahler.
%\end{enumerate}
%\end{remark}
The following result follows immediately from Theorem \ref{1uniformization}.
\begin{corollary}[cf. Corollary \ref{non-Berwald}]\label{0non-Berwald}
Let $D\subset \mathbb C^n$ be a smooth, bounded, strongly convex domain, which is not biholomorphic to $\mathbb B^n$. Its Kobayashi metric is neither complex Berwald nor weakly complex Einstein--Finsler in the sense of Aikou.
\end{corollary}
%To characterize the K\"ahler property of these  Kobayashi metrics more precisely, we introduce an intermediate condition between the weakly K\"ahler--Finsler condition and the K\"ahler--Finsler condition.
% (cf.Definition \ref{almostKF-def}), and we prove that
%\begin{theorem}[cf.Theorem \ref{almost-weaker-2}]
%In complex dimension two, a strongly pseudoconvex complex Finsler metric is weakly K\"ahler--Finsler if and only if it is
%strictly weakly K\"ahler--Finsler.
%\end{theorem}

\begin{remark}
By Corollary \ref{0non-Berwald}, for every
 smooth, bounded, strongly convex domain which is not biholomorphic to $\mathbb B^n$, its Kobayashi metric is a weakly K\"ahler--Finsler metric, which is not K\"ahler--Berwald, which, for the first time, provide smooth metrics of such kind. Such Kobayashi metrics coincide with the Carath\'eodory metrics thus have constant holomorphic sectional curvature, and they are strongly convex metrics.
\end{remark}

We also obtain the following characterization for the  locally irreducible de Rham factors of a K\"ahler--Berwald manifold $(M,F)$.
\begin{theorem}[cf. Theorem \ref{curv-sub-KB}]\label{curv-sub-KB1}
For a K\"ahler--Berwald manifold $(M,F)$, let $(U^{(j)},F^{(j)}),j=1,\cdots,k,$ be its locally irreducible de Rham factors. Then  each $(U^{(j)},F^{(j)})$ is a totally geodesic K\"ahler--Berwald sub-manifold of $(M,F)$, and for fixed  $z^{(i)},\,i\neq j,$
\begin{align}
(R^{(j)})^{\alpha_j}_{\beta_j;\gamma_j,\overline{ \sigma j}}\left(z^{(j)},v^{(j)}\right)
&=R^{\alpha_j}_{\beta_j;\gamma_j,\overline{ \sigma j}}\left(z,\pi^{\ast0}(v^{(j)})\right),\\
(R^{(j)})_{\beta_j\overline{\mu_j};\gamma_j,\overline{ \sigma j}}\left(z^{(j)},v^{(j)}\right)
&=R_{\beta_j\overline{\mu_j};\gamma_j,\overline{ \sigma j}}\left(z,\pi^{\ast0}(v^{(j)})\right).
\end{align}
Especially,
\begin{align}
K^{(j)}\left(z^{(j)},v^{(j)}\right)
&= K\left(z,\pi^{\ast0}(v^{(j)})\right),\quad \forall\, v^{(j)}\in \tilde M^{(j)}_{z^{(j)}};\\
B^{(j)}\left(z^{(j)},v^{(j)},w^{(j)}\right)
&=B\left(z,\pi^{\ast0}(v^{(j)}), \pi^{\ast0}(w^{(j)})\right),\quad \forall\, v^{(j)},w^{(j)}\in \tilde M^{(j)}_{z^{(j)}}.
\end{align}
where \[\pi^{\ast0}(v^{(j)})
=\left(0,\cdots,0,v^{(j)},0,\cdots,0\right),\,
\pi^{\ast0}(w^{(j)})
=\left(0,\cdots,0,w^{(j)},0,\cdots,0\right).\]
\end{theorem}
\begin{remark}
By Theorem \ref{curv-sub-KB}, the de Rham decomposition theorem for K\"ahler--Berwald manifolds, many results with K\"ahler geometry with holomorphic sectional/bisectional curvature conditions
can be generalized to the K\"ahler--Berwald setting, for example, the generalized Frankel conjecture verified in this paper.
\end{remark}

The rest of this subsection is about rigidity theorems for invariant metrics.
Invariant metrics, such as the Bergman metric, the Kobayashi and Carath\'eodory metrics are basic tools in the study of several complex variables, especially for bounded domain.
For example, the well known uniformization theorem of Q.-K. Lu \cite{Lu66}:
if the Bergman metric of a bounded domain $\Omega$ is complete and of constant holomorphic sectional curvature, then the $\Omega$ is biholomorphic to $\mathbb B^n$.
If consider for bounded pseudoconvex domains, there also exists uniquely the canonical complete K\"ahler--Einstein metric  \cite{Cheng-Yau,Mok-Yau}.
S.-Y. Cheng conjectured that
the Bergman metric of a smooth bounded strongly pseudoconvex domain $\Omega\in\mathbb C^n$ is K\"ahler--Einstein if and only if $\Omega$ is biholomorphic to $\mathbb B^n$ \cite{Cheng79};
and Yau conjectured that the Bergman metric of a bounded  pseudoconvex domain $\Omega$ is a complete K\"ahler--Einstein metric if and only if $\Omega$ is a bounded homogeneous domain \cite{Yau}.
Cheng's conjecture has been confirmed by X. Huang and M. Xiao in 2021 \cite{Huang-Xiao}, also previously obtained by Fu--Wong \cite{FW97} and Nemirovski--Shafikov \cite{NS06} in the case of complex dimension two.

There are also rigidity theorems for the Kobayashi and Carath\'eodory metrics.
In 1977, B. Wong proved the following theorem.
\begin{theorem}[\cite{Wong77}]\label{Wong77}
Let $M$ be a simply-connected complete hyperbolic complex manifold with $F_K =F_C$, assumed $F_K$ to be a K\"ahler metric. Then, $M$ is biholomorphic to $\mathbb B^n$.
\end{theorem}
Stanton subsequently proved the following more general result in 1983.
\begin{theorem} [\cite{Stanton}]\label{Stanton}
 Let $M$ be a complete hyperbolic complex manifold. Assume that there is a point on $M$ at which $F_K$ and $F_C$ are equal, and one of these metrics is Hermitian and of class $C^\infty$. Then
$M$ is biholomorphic to $\mathbb B^n$.
\end{theorem}

In this paper, we obtain the following rigidity theorem.
\begin{theorem}[cf. Theorem \ref{Kob-KB}]\label{0Kob-KB}
 For a complex manifold $M$, we have
\begin{itemize}
  \item [1)] $F_K$ is K\"ahler--Berwald if and only if it is K\"ahler;
  \item [2)] $F_C$ is K\"ahler--Berwald if and only if it is K\"ahler, when $M$ is in additionally assumed simply connected.
\end{itemize}
\end{theorem}

By Theorems \ref{Stanton}--\ref{0Kob-KB}, we obtain the following result.
\begin{theorem}\label{glob-rig}
 Let $M$ be a complete hyperbolic complex manifold. Assume that there is a point on $M$ at which $F_K$ and $F_C$ are equal,
then the following are equivalent
 \begin{itemize}
   \item [1)] $F_K$ is K\"ahler--Berwald and of class $C^\infty$;
   \item [2)] $F_C$ is K\"ahler--Berwald and of class $C^\infty$, and $M$ is simply connected;
   \item [3)] $M$ is biholomorphic to $\mathbb B^n$.
 \end{itemize}
\end{theorem}
\begin{remark}
There are many bounded domain such that $F_K=F_C$ at some points, especially for bounded convex domains and strongly linear convex domains \cite{Lem81,Lem82,Kos-War13}.
\end{remark}

However, there exists smooth bounded strongly pseudoconvex domain such that $F_K\neq F_C$, for example \cite{Burns-Schnider}:
\[\mathcal S:=\{(z,w)\in \mathbb C^2: e^{-\pi}<\vert z\vert<1, \, \sin(\log\vert z\vert)+\vert w\vert^2<0\}.\]
Since $\mathcal S$ has spherical boundary, thus it is universally covered by $\mathbb B^2$ \cite{Yu}, hence the Kobayashi metric of $\mathcal S$ is K\"ahler. By Theorem \ref{Stanton}, $F_K\neq F_C$ at any point, since $\mathcal S$ is not simply connected, thus do not biholomorphic to $\mathbb B^2$.

For general smooth bounded strongly pseudoconvex domain, the two invariant metrics agree in a nonempty open set of tangent directions at points sufficiently close to the boundary \cite{BFW19,Gau-Zim,BK94}.
The following Theorem~\ref{K=C} was established under $C^3$-smooth boundary condition by Bracci--Forn\ae ss--Wold \cite{BFW19}, and for $C^2$ case by Gaussier--Zimmer \cite{Gau-Zim}. A version for $C^6$ case also follows from the arguments in Section~4 of \cite{BK94}.

\begin{theorem}[\cite{BK94,BFW19,Gau-Zim}]\label{K=C}
Suppose that $\Omega\subset \mathbb C^d$ is a bounded strongly pseudoconvex domain with $C^2$ boundary.
There exists $\varepsilon>0$ such that, whenever
$\delta_\Omega(z):=\operatorname{dist}(z,\partial\Omega)<\varepsilon$, there is a nonempty open set
$E_\Omega(z)\subset T_z^{1,0}\Omega\setminus\{0\}$ on which the
Kobayashi and Carath\'eodory metrics agree.
\end{theorem}

Based on Theorem \ref{K=C}, Gaussier--Zimmer obtained the following rigidity theorem using the K\"ahler condition \cite{Gau-Zim}.
\begin{theorem} [\cite{Gau-Zim}]\label{Gaussier-Zimmer}
%Suppose that $M$ is a Stein manifold with $dim_{\mathbb C}M=n\geq2$ and $g$ is a complete K\"ahler metric on $M$. If there exists a compact set $\Sigma\subset M$ such that $M\setminus\Sigma$ is connected and $g$ has constant negative holomorphic sectional curvature on  $M\setminus\Sigma$, then the universal cover of $M$ is biholomorphic to $\mathbb B^n$.
Let $\Omega\subset \mathbb C^n$ be a bounded strongly pseudoconvex domain with $C^2$-boundary. Then the following are equivalent:
\begin{itemize}
  \item [1)] the Kobayashi metric of $\Omega$ is a K\"ahler metric on $U \cap \Omega$, where $U\supseteq \partial \Omega$ is an open set;
  \item [2)] the Kobayashi metric of $\Omega$ is a K\"ahler metric with constant holomorphic sectional curvature;
  \item [3)] the universal cover of $\Omega$ is biholomorphic to $\mathbb B^n$.
\end{itemize}
\end{theorem}

Similar to the Kobayashi case, R. X. Dong, R. Y. Wang and B. Wong obtained the following rigidity theorem for the Carath\'eodory metric \cite{Dong-Wang-Wong}.
\begin{theorem}[\cite{Dong-Wang-Wong}]\label{DWW}
 Suppose that $\Omega \subset \mathbb{C}^{n},$ $n \geq2$, is a bounded strictly pseudoconvex domain with a $C^{2}$-boundary. Then, $\Omega$ is biholomorphic to $\mathbb B^n$ if and only if
 the Carath\'eodory metric of $\Omega$ is a K\"ahler metric on $U \cap \Omega$, where $U\supseteq \partial D$ is an open set.
\end{theorem}

By Theorem \ref{0Kob-KB}, and the de Rham decomposition theorem for K\"ahler--Berwald manifolds (Theorem \ref{1characKB}),
H. Gaussier and A. Zimmer's result can be generalized to the K\"ahler--Berwald setting.
\begin{corollary}[cf. Corollary \ref{KB-Ball}]\label{1Kob-Ball}
Let $\Omega\subset \mathbb C^n,n\geq2,$ be a bounded strongly pseudoconvex domain with $C^2$-boundary, and assume that its Kobayashi metric $F_K$ is strongly pseudoconvex. Then the following are equivalent:
\begin{enumerate}
    \item [1)] the universal cover of $\Omega$ is biholomorphic to $\mathbb B^n$;
  \item  [2)] $F_K$ is K\"ahler--Einstein and of constant holomorphic sectional curvature;
  \item  [3)] $F_K$ is weakly K\"ahler--Finsler, weakly complex Einstein--Finsler on $U \cap \Omega$, where $U\supseteq \partial \Omega$ is an open set;
  \item  [4)] $F_K$ is K\"ahler--Berwald on $U \cap \Omega$, where $U\supseteq \partial \Omega$ is an open set.
\end{enumerate}
\end{corollary}
The Kobayashi metrics of bounded strongly convex domains are weakly K\"ahler--Finsler (cf. \cite{Lem81,Lem82,Pat,AP2}). Thus, by Corollary \ref{1Kob-Ball},
\begin{corollary}\label{KFE-Ball2}
Let $\Omega\subset \mathbb C^n,n\geq2,$ be a smooth bounded strongly convex domain. Then the following are equivalent:
\begin{enumerate}
  \item  [1)] $\Omega$ is biholomorphic to $\mathbb B^n$;
  \item  [2)] $F_K$ is weakly complex Einstein--Finsler;
  \item  [3)] $F_K$ is complex Berwald.
\end{enumerate}
\end{corollary}
Similarly, by Theorem \ref{0Kob-KB}, and Theorem \ref{1characKB}, Theorem \ref{DWW} can be generalized to the K\"ahler--Berwald setting.
\begin{corollary}[cf. Corollary \ref{DWW-KB}]\label{1DWW-KB}
Let $\Omega \subset \mathbb{C}^{n},$ $n \geq2$, be a simply connected, bounded strongly pseudoconvex domain with $C^{2}$-boundary. Then the following are equivalent:
\begin{enumerate}
  \item [1)]  $\Omega$ is biholomorphic to $\mathbb B^n$;
  \item [2)]  The Carath\'eodory metric of $\Omega$ is a strongly pseudoconvex K\"ahler--Berwald metric on $U \cap \Omega$, where $U\supseteq \partial \Omega$ is an open set;
  \item  [3)] The Carath\'eodory metric of $\Omega$ is a strongly pseudoconvex weakly K\"ahler--Finsler, weakly complex Einstein--Finsler metric on $U \cap \Omega$, where $U\supseteq \partial \Omega$ is an open set.
\end{enumerate}
\end{corollary}

On convex domain (bounded or unbounded), the Kobayashi metric and Carath\'eodory metrics coincide, thus a corollary of Theorem \ref{0Kob-KB} is as follows.
\begin{corollary}[cf. Corollary \ref{rig-unb-con-Kob}]
For every convex domain $\Omega\subset \mathbb C^n$, its Kobayashi metric is K\"ahler--Berwald if and only if $\Omega$ is biholomorphic to $\mathbb B^n$.
\end{corollary}

The rest of this paper is organized as follows:
in Section \ref{Sec-cpx-F}, we give a brief introduction to complex Finsler geometry;
in Section \ref{sec:KF=KB}, we introduce the formally stronger yet equivalent weakly K\"ahler--Finsler condition (SWKF), and
prove Theorem \ref{thm:KF=KB1}, the equivalence between the K\"ahler--Finsler and K\"ahler--Berwald conditions;
in Section \ref{sec-loc-symm}, we prove Theorem \ref{thm}, Theorem \ref{1pal-eq}, and results about locally symmetric complex Finsler metrics;
in Section \ref{sec-de-Rham}, we prove Theorems \ref{1characKB}--\ref{KBFrankel}, the de Rham decomposition theorem, and characterization of K\"ahler--Berwald manifold, and verify the generalized Frankel conjecture for the K\"ahler--Berwald setting;
in Section \ref{sec-unif}, we prove Theorem \ref{1uniformization}, the uniformization theorem, and Theorem \ref{0Kob-KB}, the K\"ahler--Berwald rigidity for the Kobayashi metric, and its corollaries.

\section{Complex Finsler geometry}\label{Sec-cpx-F}
	
In this section, we briefly review some notions and results for strongly pseudoconvex complex Finsler manifold $(M,F)$.

We always use the Einstein convention in this paper.
And we adopt the convention of using a subscript to indicate differentiation of $G$ and a semicolon to distinguish derivatives with respect to the base coordinates from those with respect to the fiber coordinates. For instance,
$G_{j\bar{k};i\bar l} := \partial_i\partial_{\bar l}\dot\partial_{k}\dot\partial_{\bar k}G$,
where
\[
\partial_i:=\frac{\partial}{\partial z^i},\,
\partial_{\bar l}:=\frac{\partial}{\partial \bar z^l},\,
\dot\partial_{ j}:=\frac{\partial}{\partial v^j},\,
\dot\partial_{\bar k}:=\frac{\partial}{\partial \bar{v}^k}.
\]

Denote $(G^{\bar{j}l}):=(G_{i\bar{j}})^{-1}$,
%the inverse matrix of $(G_{i\bar{j}})$ by $(G^{\bar{j}l})$,
and $\mathcal V:=span_{\mathbb C}\{\dot\partial_1,\cdots,\dot\partial_n\}=ker (d\pi)\subset T^{1,0}\tilde M $ for the vertical sub-bundle over $\tilde M$, where $\pi: \tilde M \rightarrow M$ is the canonical projection. Then $F$ induces a Hermitian metric on $\mathcal V:$
\begin{equation}
\langle V,W\rangle:=G_{a\bar b}(z,v)V^a\bar W^b,\quad\forall \,V,W\in \mathcal V_{(z,v)}.\label{vert-metric}
\end{equation}

The Chern--Finsler connection of $F$ is defined as follows.

\begin{definition}[\cite{AP}]\label{Abate-Chern}
Let $F$ be a strongly pseudoconvex complex Finsler metric of $M$, and $\langle \cdot,\cdot\rangle$ is defined as \eqref{vert-metric}. The Chern--Finsler connection is the unique complex vertical connection
$\nabla:\Gamma(\mathcal V)\rightarrow\Gamma(T^\ast_{\mathbb C}\tilde M\otimes \mathcal V ) $ such that
\begin{equation}
X\langle V,W\rangle=\langle\nabla_X V,W\rangle+\langle V,\nabla_{\bar X} W\rangle \label{Chern--Finsler}
\end{equation}
for all $X\in T^{1,0}\tilde M$ and $V,W\in \Gamma(\mathcal V).$
Furthermore, this connection is good.
\end{definition}

In local coordinates, the connection form of $\nabla$ is defined by
\[
\omega^i_j=G^{\bar c i}\bar\partial G_{j\bar c}
=:\Gamma^i_{j;k}dz^k+\Gamma^i_{jk}\psi^k,
\]
where
\begin{equation}
\Gamma_{j;k}^{i} = \dot\partial_j \Gamma_{;k}^{i},\,
\Gamma_{;k}^{i}(z,v) =  G^{\bar{s}i} G_{\bar{s};k},\,
\Gamma_{jk}^{i} = G^{\bar{s}i} G_{j\bar{s}k}.
\label{con-coe}
\end{equation}
And $\Gamma_{j;k}^{i}, \Gamma_{;k}^{i}$ and $\Gamma_{jk}^{i}$
are called the horizontal, nonlinear and vertical connextion coefficients of the Chern--Finsler connection of $F$, respectively.
Denoting
\[\psi^i:=dv^i+\Gamma^i_{;j}dz^j,\quad \delta_i:=\partial_i-\Gamma^j_{;i}\dot\partial_j, \]
then $\{dz^i,\psi^j\}_{i,j=1}^n$ is the locally dual coframe of $\{\delta_i,\dot\partial_j\}_{i,j=1}^n.$
And $\{\delta_1,\cdots,\delta_n\}$ forms a local frame in $\mathcal H$, a horizontal sub-bundle of $T^{1,0}\tilde M$ such that $T^{1,0}\tilde M=\mathcal H\oplus\mathcal V.$
Using the homogeneous property of $G(z,v)$ with respect to $v$, by the Euler homogeneous identity, it is easy to check (cf. \cite{AP})
\begin{eqnarray}
G(z,v)&=&G_{a}(z,v)v^a=G_{\bar a}(z,v)\bar v^a;\\
0&=&G_{ab}(z,v)v^a=G_{ab\bar c }(z,v)v^a; \\ \Gamma^i_{j;k}(z,v)v^j&=&\Gamma^i_{;k}(z,v),\,
\Gamma^i_{jk}(z,v)v^j=0.\label{AP-homo}
\end{eqnarray}

Using the bundle isomorphism $\Theta:\mathcal V\rightarrow\mathcal H,\, \frac{\partial}{\partial v^j}\mapsto\delta_j$, $\nabla$ induces a connection for the horizontal sub-bundle:
\[\nabla_X H:=\Theta[\nabla_X(\Theta^{-1}(H))],\quad \forall \,H\in \Gamma(\mathcal H),\forall \,X\in T\tilde M.\]
Thus we can also view the Chern--Finsler connection as defined on the horizontal sub-bundle.

By Corollary 2.4.2 in \cite{AP}, for a weakly K\"ahler--Finsler manifold $(M,F)$, a curve $\gamma$ is a geodesic if it satisfies
\begin{equation}
\nabla_{T^{\mathcal{H}}+\bar T^{\mathcal{H}}} T^{\mathcal{H}} = 0,
\label{geodesic}
\end{equation}
where $\dot\gamma(t) := \dfrac{\partial \gamma}{\partial t}(t)$, and $T^{\mathcal{H}} := \chi(\dot\gamma)=\dot\gamma^i\delta_i$ is the horizontal lifting of $\dot\gamma$, from the tangent bundle $T^{1,0}M$ to the horizontal sub-bundle $\mathcal{H}$. Locally, equation \eqref{geodesic} can be written as \cite{AP}
\begin{equation}
\ddot\gamma^{i} + \Gamma_{j;k}^{i}(\gamma,\dot\gamma) \dot\gamma^{j} \dot\gamma^{k} = 0, \quad i,j,k = 1,\dots,n.
\end{equation}
This is also the geodesic equations for K\"ahler metrics when %viewing $ \Gamma_{j;k}^{i}(\gamma;\dot\gamma)$ as Chern connection coefficients when
$\Gamma_{j;k}^{i}(\gamma,\dot\gamma)=\Gamma_{j;k}^{i}(\gamma)$. Similarly, the parallel vector field equation along $\gamma$ is given by
\begin{equation}
\dot\xi^{i} + \Gamma_{j;k}^{i}(\gamma,\dot\gamma) \dot\gamma^{j} \xi^{k}= 0, \qquad i,j,k = 1,\dots,n.
\label{geod}
\end{equation}

The curvature operator of the Chern--Finsler connection of $F$ is defined by $\Omega(X,Y)Z:=R_X(X,Y).$
Locally, $\Omega =\Omega_j^i\otimes[dz^j\otimes\delta_i+\psi^j\otimes \dot\partial_i]$
with
\[
\Omega_j^i=\bar\partial \omega_j^i
=R_{j;k\bar l}^i d z^k\wedge d\bar z^l
+R_{jk;\bar l}^i \psi^k\wedge d\bar z^l
+R_{j\bar k;l}^i d z^l\wedge \overline{ \psi^k}
+R_{jk\bar l}^i \psi^k\wedge \overline{\psi^l}.
\]
The $h\bar h$-part curvature is
\begin{equation}
R_{j;k\bar l}^i
:=-\delta_{\bar l}\Gamma^i_{j;k}-\Gamma^i_{js}\delta_{\bar l }\Gamma^s_{;k},\label{HH-curv}
\end{equation}
also one can denote
\begin{equation}
R_{j \bar s;k\bar l}:=G_{i\bar s}R_{j;k\bar l}^i.\label{HH-curv'}
\end{equation}
The $h\bar v$-part curvature is
\begin{equation}
R_{j\bar l;k}^i
:=-\dot\partial_{\bar l}\Gamma^i_{j;k}-\Gamma^i_{js}\dot\partial_{\bar l}\Gamma^s_{;k}.\label{HV-curv}
\end{equation}

When $F$ comes from a Hermitian metric, the non-vanishing part curvature is only the $h\bar h$-part, which then reduce to the Chern curvature.
Thus the holomorphic sectional curvature $K$ for a general complex Finsler metric $F$ along direction $(z,v)\in\tilde{M}$ is defined by \cite{AP}
\begin{equation}
K(z,v)=\frac{2\langle \Omega(\chi,\overline{\chi})\chi,
\chi\rangle_{(z,v)}}{G^2(z,v)}
=-\frac{2}{G^2(z,v)}\delta_{\bar l}\Gamma^i_{j;k}v^jv^k\bar v^l
,\label{hsc-a}
\end{equation}
where $\chi=\sum_{i=1}^nv^i\delta_i$ is the complex horizontal radial vector field.
	
The holomorphic bisectional curvature of $F$ with respect to $(z,v),(z,w)\in \tilde{M}$ is defined as \cite{SZ}	
\begin{equation}
B(z,v,w)=\frac{2\left\langle \Omega\left(\chi_v(w),\overline{\chi_v(w)}\right)
\chi,\chi\right\rangle_{(z,v)}}{G(z,v)G(z,w)},\label{hbsc-a}
\end{equation}
where $\chi_v(w)=\sum_{i=1}^nw^i\delta_i\vert_{(z,v)}$ is the horizontal lifting of $w$ along $v$.

Let  $(E, F)$ be a complex Finsler bundle over a Hermitian manifold $(M, g)$, where $g=\sum_{\alpha, \beta} g_{\alpha \bar{\beta}}(z) d z^{\alpha} \otimes d \bar{z}^{\beta}$. The mean curvature $K_{i \bar{j}}$ of $(E, F)$ over $(M, g)$ was introduced by Aikou in  \cite{Ai97} as
\[K_{j}^{i}:=\sum_{\alpha, \overline{\beta}} g^{\overline{\beta} \alpha} R_{j; \alpha \overline{\beta}}^{i}.\]
$(E,F)$ is called weakly complex Einstein--Finsler if $K_{j}^i=\varphi(z)\delta^i_j,$
and it is called to be complex Einstein--Finsler if $K_{j}^i=c\delta^i_j$.
In \cite{SZ18}, L. Sun and C. Zhong proved the following result. \begin{theorem}[\cite{SZ18}]\label{SZ18EF}
Let  $F$ be a strongly pseudoconvex complex Finsler metric. If $F$ is a weakly complex Einstein--Finsler metric in the sense of Aikou, it is a complex Berwald metric.
\end{theorem}
Thus by \eqref{XZ-B-wK-K} and Theorem \ref{SZ18EF}, the following corollary holds.
\begin{corollary}[\cite{SZ18}]\label{SXZ}
Let  $F$ be a strongly pseudoconvex weakly K\"ahler--Finsler metric. If $F$ is weakly complex Einstein--Finsler in the sense of Aikou, it is K\"ahler--Berwald.
\end{corollary}

At last, in order to generalize the canonical property for symmetric Hermitian metrics into Finsler circumstance in some generalized sense, we recall the Chern--Finsler nonlinear connection.

%\begin{definition}\label{Chern--Finsler nonlinear-class}
%Two strongly pseudoconvex complex Finsler metrics are called in a Chern--Finsler (or Chern--Finsler nonlinear, since the complex Rund connection coincides with the horizontal Chern--Finsler connection \cite{Rund}) class, if they share same horizontal Chern--Finsler connection.
%\end{definition}

\begin{definition}\label{nonl-con}
For a strongly pseudoconvex complex Finsler manifold $(M,F)$,
the Chern--Finsler nonlinear connection is the map $\nabla:\Gamma(T^{(1,0)}M)\rightarrow \Gamma (T^\ast_{\mathbb C}M\otimes T^{(1,0)}M)$ defined by:
\[
\nabla\xi:=\left(d\xi^i+\xi^j\Gamma^i_{j;k}(\xi)dz^k  \right)\frac{\partial}{\partial z^i},\, \xi\in \Gamma(\tilde M),
\]
where $\Gamma^i_{j;k}(\xi)$ is the horizontal part of the Chern--Finsler connection of $F$.
\end{definition}
\begin{remark}
For a strongly pseudoconvex complex Berwald metric, $\Gamma^i_{j;k}(z,v)=\Gamma^i_{j;k}(z),\forall \,0\neq v\in T^{(1,0)}_z M,$
$\nabla$ is indeed a linear connection. We just still call $\nabla$ the Chern--Finsler nonlinear connection of $F$, and it coincides with the Chern connection of a Hermitian metric in this case.
\end{remark}

\begin{remark} Notice that two irreducible K\"ahler metrics with same Chern connection differ only by a multiplicative constant. Thus this definition is trivial for irreducible K\"ahler metrics. That is, on an irreducible K\"ahler--Berwald manifold $(M,F)$, there is only one K\"ahler metric, up to multiplicative constants, whose Chern connection coincides with the Chern--Finsler nonlinear connection of $F$.
In this case, the number of K\"ahler--Berwald metrics on $M$ with different Chern--Finsler nonlinear connection coincides with the number of K\"ahler metrics with different Chern connection. In particular, irreducible locally symmetric Hermitian metric is unique up to scaling.
 \end{remark}

 And by definition, one has the following results.
\begin{lemma}\label{Chern--Finsler nonlinear-geod}
On a K\"ahler--Berwald manifold, for two K\"ahler--Berwald metrics share same Chern--Finsler nonlinear connection,
the geodesics, parallel vector fields, the holonomy groups, and the horizontal part of the Chern--Finsler curvature operator of such two metrics are coincide.
\end{lemma}

\section{Generalizations of the K\"ahler condition in complex Finsler geometry}\label{sec:KF=KB}

In this section, we consider the relationships between weakly K\"ahler--Finsler, the formally stronger condition than weakly K\"ahler--Finsler, K\"ahler--Finsler, strongly K\"ahler--Finsler and K\"ahler--Berwald.
Below we denote the above conditions, respectively, as
\[WKF,\,SWKF,\,KF,\, SKF\,\text{and} \,KB.\]

Recall the $WKF, SWKF, KF, \,\text{and}\, SKF$ conditions are defined, respectively, as
 \begin{align}
G_a(z,v)\left(\Gamma^a_{b;c}(z,v)-\Gamma^a_{c;b}(z,v)\right)v^b&=0,
\quad\forall \,0\neq v\in T_z^{1,0}M;\label{weakly-KF}\\
G_a(z,v)\left(\Gamma^a_{b;c}(z,v)-\Gamma^a_{c;b}(z,v)\right)&=0,
\quad\forall \,0\neq v\in T_z^{1,0}M;\label{new-KF}\\
\left(\Gamma^a_{b;c}(z,v)-\Gamma^a_{c;b}(z,v)\right)v^b&=0,
\quad\forall \,0\neq v\in T_z^{1,0}M;\\
\Gamma^a_{b;c}(z,v)-\Gamma^a_{c;b}(z,v)&=0,
\quad\forall \,0\neq v\in T_z^{1,0}M. \label{Recall-KF}
\end{align}
By the Euler homogeneous identity, $G_a=G_{a\bar l}\bar v^l,$ thus
$G_aG^{\bar ka}=\bar v^k,$ and
\[
G_a\Gamma^a_{b;c}=G_aG^{\bar ka}\delta_cG_{b\bar k}
=\bar v^{k}\delta_cG_{b\bar k}
=\delta_cG_{b}.
\]
Hence $WKF, SWKF, KF \, \text{and}\, SKF$ also can be written, respectively, as
\begin{align}
(\delta_cG_{b}-\delta_bG_{c})v^c&=0,\quad\forall \,0\neq v\in T_z^{1,0}M;\\
\delta_cG_b-\delta_bG_c&=0,\quad
\forall \,0\neq v\in T_z^{1,0}M;\label{new-KF'}\\
(\delta_cG_{b\bar a}-\delta_bG_{c\bar a})v^c&=0,\quad\forall \,0\neq v\in T_z^{1,0}M;\\
 \delta_cG_{b\bar a}-\delta_bG_{c\bar a}&=0,\quad\forall \,0\neq v\in T_z^{1,0}M.
 \end{align}
To view the relationships between $SWKF$ and $KF (= SKF)$ more clearly, \eqref{new-KF'} can also be written as
 \begin{eqnarray}
 (\delta_cG_{b\bar a}-\delta_bG_{c\bar a})\bar v^a&=&0,\quad\forall \, 0\neq v\in T_z^{1,0}M.\label{new-KF2}
 \end{eqnarray}

A nature question is as follows.
\begin{question}\label{almost-KF}
Does $WKF\Leftrightarrow SWKF$ or $SKF\Leftrightarrow SWKF$
hold?
\end{question}
We confirm the first part of Question \ref{almost-KF} in case of complex dimension two.
\begin{proposition}\label{new-KF-th}
Let $M$ ba a complex manifold of complex dimension two. Then every strongly pseudoconvex weakly K\"ahler--Finsler metric on $M$ satisfies \eqref{new-KF}.
\end{proposition}
\begin{proof}
By condition, we have
\begin{equation}
G_a\left(\Gamma^a_{b;c}-\Gamma^a_{c;b}\right)v^b=0,
\quad\forall \,0\neq v\in T_z^{1,0}M.\label{neq-KF-eq1}
\end{equation}
In case of complex dimension-$2$, \eqref{neq-KF-eq1} is as follows
\begin{equation}
G_a\left[\left(\Gamma^a_{1;c}-\Gamma^a_{c;1}\right)v^1
+\left(\Gamma^a_{2;c}-\Gamma^a_{c;2}\right)v^2\right]
=0.\label{neq-KF-eq2}
\end{equation}
For $c=1$, \eqref{neq-KF-eq2} reduces to
\begin{eqnarray}
G_a\left[\left(\Gamma^a_{1;1}-\Gamma^a_{1;1}\right)v^1
+\left(\Gamma^a_{2;1}-\Gamma^a_{1;2}\right)v^2\right]
=0,\quad\forall \,0\neq v\in T_z^{1,0}M.\label{neq-KF-eq3}
\end{eqnarray}
For $c=2$, \eqref{neq-KF-eq2} reduces to
\begin{eqnarray}
G_a\left[\left(\Gamma^a_{1;2}-\Gamma^a_{2;1}\right)v^1
+\left(\Gamma^a_{2;2}-\Gamma^a_{2;2}\right)v^2\right]
=0,\quad\forall \,0\neq v\in T_z^{1,0}M.\label{neq-KF-eq4}
\end{eqnarray}
\eqref{neq-KF-eq3} and \eqref{neq-KF-eq4} provide
\begin{eqnarray}
G_a\left(\Gamma^a_{2;1}-\Gamma^a_{1;2}\right)v^2
=G_a\left(\Gamma^a_{1;2}-\Gamma^a_{2;1}\right)v^1=0,
\quad\forall \,0\neq v\in T_z^{1,0}M,\label{neq-KF-eq5}
\end{eqnarray}
Notive that for any $v=(v^1,v^2)\neq0,$ we have $v^1\neq0$ or $v^2\neq0.$
Thus, by \eqref{neq-KF-eq5}, we have
\begin{eqnarray}
G_a\left(\Gamma^a_{1;2}-\Gamma^a_{2;1}\right)
=0.\label{neq-KF-eq7}
\end{eqnarray}
And obviously, \eqref{neq-KF-eq7} gives the following equation for all indices b and c:
\[G_a\left(\Gamma^a_{b;c}-\Gamma^a_{c;b}\right)
=0.
\]
\end{proof}

%Especially, since the Kobayashi metric of a smooth bounded strongly convex domain is weakly K\"ahler--Finsler, by Proposition \ref{new-KF-th}, we have
%\begin{corollary}\label{new-KF-Kob}
%The Kobayashi metric of any smooth bounded strongly convex domain in $\mathbb C^2$ satisfies \eqref{new-KF}.
%\end{corollary}
%By Theorem \ref{non-Berwald} and Corollary \ref{new-KF-Kob},
\begin{remark}
By Theorem \ref{new-KF-th}, the uniformization theorem does not hold in the $WSKF$ setting.
The answer to  Question \ref{almost-KF} is know clear.
\begin{itemize}
  \item [1)] The coincidence $WKF\Leftrightarrow SWKF$ for general dimensional cases has been confirmed by the first-named author in \cite{Ge-wKF}, after we obtained the present Theorem \ref{new-KF-th};
  \item [2)] The relationship $SWKF\nRightarrow KF$ is implied by Theorem \ref{new-KF-th} together with the coincidence $KF\Leftrightarrow KB$, established below in Theorem \ref{thm:KF=KB}.
\end{itemize}
\end{remark}

Now, we recall relationships between $WKF(=SWKF)$ and $KF(=SKF).$
In \cite{CS}, B. Chen and Y. Shen proved the coincidence $KF\Leftrightarrow SKF$. %, they also obtained a condition for weakly  K\"ahler--Finsler metric to be K\"ahler--Finsler.
%\begin{theorem}[{\cite[Theorem 4.1]{CS}}]\label{CS-van}
%If the weak $h\bar v$-curvature of the Chern--Finsler connection of a complex Finsler metric $F$ vanishes identically, then $F$ is K\"ahler if and only if it is weakly K\"ahler.
%\end{theorem}
%Recall that the $h\bar v$-part curvature defined by \eqref{HV-curv} is as follows:
%\begin{equation}
%R_{j\bar l;k}^i
%:=-\dot\partial_{\bar l}\Gamma^i_{j;k}-\Gamma^i_{js}\dot\partial_{\bar l}\Gamma^s_{;k},\label{CS-hv-curv}
%\end{equation}
%The $h\bar v$-part curvature vanishing condition leads to
%\begin{equation}
%0=R_{j\bar l;k}^iv^jv^k
%=-\dot\partial_{\bar l}\Gamma^i_{j;k}v^jv^k,\label{hv-B-cond}
%\end{equation}
%where in the first equality we used $\Gamma^i_{js}v^j=G^{\bar ki}G_{\bar kjs}v^j=0$, applying the Euler formula to the $(0,0)$-homogeneous map $G_{\bar ks}$.

%Thus \eqref{hv-B-cond} implies $\mathbb G^i_{\bar j}:=\dot\partial_{\bar j}\mathbb G^i=0,$
%a holomorphic condition for $\mathbb G^i$.
%In \cite{Ge-wKF}, the first-named author obtained the coincidence $WKF\Leftrightarrow SWKF.$

%Based on the new condition $SWKF$, he obtained the following relationship between weakly K\"ahler--Finsler and K\"ahler--Finsler, however, this .
In \cite{XZ}, H. Xia and C. Zhong obtained the following result.
\begin{lemma}[{\cite[Theorem 1.1]{XZ}}]\label{lem}
A weakly K\"ahler--Finsler metric $F$ is K\"ahler--Finsler if and only if
\[G_{ij}\mathbb G^i_{\bar l}=0.\]
\end{lemma}
%For completeness, we provide a short proof of Lemma \ref{lem}.
%\begin{proof}
%By the Euler formula for homogeneous maps and direct calculation,
%\begin{eqnarray*}
%G_i\dot\partial_{\bar l}\Gamma^i_{;k}
%&=&G_i\dot\partial_{\bar l}\left(G^{\bar ai}G_{\bar a;k}\right)\\
%&=&G_i\left( -G^{\bar ab}G_{b\bar c\bar l}G^{\bar ci}G_{\bar a;k}+G^{\bar ai}G_{\bar a\bar l;k}\right)\\
%&=&-G^{\bar ab}G_{b\bar c\bar l}v^{\bar c}G_{\bar a;k}+v^{\bar a}G_{\bar a\bar l;k}\\
%&=&0.
%\end{eqnarray*}
%Thus
%\begin{equation}
%0=\dot\partial_j\left( G_i\dot\partial_{\bar l}\Gamma^i_{;k}\right)
%=G_{ij}\dot\partial_{\bar l}\Gamma^i_{;k}
%+G_i\dot\partial_{\bar l}\Gamma^i_{j;k}.\label{lem-van1}
%\end{equation}
%By Chen--Shen's proof of Theorem \ref{CS-van},
%weakly K\"ahler--Finsler condition provides  \cite[equation $(4.2)$]{CS}:
%\begin{equation}
%G_{i\bar l}\left(\Gamma^i_{j;k}-\Gamma^i_{k;j}\right)v^k
%=-G_i\dot\partial_{\bar l}\Gamma^i_{j;k}v^k.\label{lem-van2}
%\end{equation}
%\eqref{lem-van1} and \eqref{lem-van2} imply
%\begin{equation}
%G_{i\bar l}\left(\Gamma^i_{j;k}-\Gamma^i_{k;j}\right)v^k
%=G_{ij}\dot\partial_{\bar l}\Gamma^i_{;k}v^k
%=G_{ij}\mathbb G^i_{\bar l}.\label{lem-van3}
%\end{equation}
%Thus
%\[G_{ij}\mathbb G^i_{\bar l}=0
%\Leftrightarrow G_{i\bar l}\left(\Gamma^i_{j;k}-\Gamma^i_{k;j}\right)v^k
%=0 \Leftrightarrow \left(\Gamma^i_{j;k}-\Gamma^i_{k;j}\right)v^k
%=0,\]
%since $G_{i\bar l}$ is non-vanishing.
%And $\left(\Gamma^i_{j;k}-\Gamma^i_{k;j}\right)v^k
%=0$ is just the K\"ahler--Finsler condition.
%\end{proof}
Were $\mathbb G^i_{\bar l}:=\dot\partial_{\bar l}\mathbb G^i$, and $\mathbb G^i$ is coefficient of the complex geodesic spray defined by
\[
\mathbb G^i:=\Gamma^i_{;j}v^j,\, \Gamma^i_{;j}:=G^{\bar ai}G_{\bar a;j}.
\]

Lemma \ref{lem} was used in providing weakly K\"ahler--Finsler examples that are not K\"ahler--Finsler, and K\"ahler--Finsler rigidity theorems for the Kobayashi metric \cite{Ge-wKF,Ge-KF,LLQ26}.
In the following, we prove that the K\"ahler condition implies the $\mathbb G^i_{\bar l}=0$, then we obtain the K\"ahler--Berwald conclusion.
\begin{theorem}\label{thm:KF=KB}
K\"ahler--Finsler metrics are K\"ahler--Berwald.
\end{theorem}
\begin{proof}
%For every strongly pseudoconvex complex Finsler metric,
%\begin{equation}
%\begin{aligned}
%G_i\mathbb G^i_{\bar j}
%&=G_i\dot\partial_{\bar j}\left(G^{\bar\tau i}G_{\bar \tau;l}\right)v^l\\
%&=G_i\left(-G^{\bar\tau a}G_{a\bar b\bar j}G^{\bar b i}G_{\bar \tau;l}+G^{\bar\tau i}G_{\bar j\bar \tau;l}\right)v^l\\
%&=\left(-G^{\bar\tau a}G_{a\bar b\bar j}\overline{v^{ b} }G_{\bar \tau;l}+\overline{v^{\tau} }G_{\bar j\bar \tau;l}\right)v^l\\
%&=0,
%\end{aligned}\label{Van0}
%\end{equation}
%where in the last equality we using the Euler formula for homogeneous maps. Thus
%\[0=\dot\partial_{\bar k}\left(G_{i}\mathbb G^i_{\bar j}\right)
%=G_{i\bar k}\mathbb G^i_{\bar j}+G_{i}\dot\partial_{\bar k}\dot\partial_{\bar j}\mathbb G^i,
%\]
%and
%\[
%G_{i\bar k}\mathbb G^i_{\bar j}=-G_{i}\dot\partial_{\bar k}\dot\partial_{\bar j}\mathbb G^i
%=G_{i\bar j}\mathbb G^i_{\bar k}.
%\]
%The last equality follows by interchanging $j,k$ for derivatives of $\mathbb G^i$.

Take
\[\bar{\chi}:=\overline{ v^k\delta_{ k}}
=\overline{v^k\left(\partial_k
-\Gamma_{;k}^l\dot\partial_l\right)}
=\bar v^k\partial_{\bar k}
-\overline{\mathbb G^k}\dot\partial_{\bar k}
,\]
 then
\begin{equation}
\begin{aligned}
 \bar{\chi} G_i
=\bar v^k\left(G_{i;\bar k}
-G^{ \bar l\tau}G_{\tau;\bar k}G_{i\bar l}\right)
=\bar v^k\left(G_{i;\bar k}-G_{i;\bar k}\right)
=0.
\end{aligned}\label{frak-D}
\end{equation}
 Differentiated in $v^j$ for \eqref{frak-D}, yields
\[ 0=\dot\partial_j\left(\bar{\chi} G_i\right)
=\bar{\chi}\dot\partial_j G_i
-\dot\partial_j\overline{\mathbb G^k}G_{i\bar k}
=\bar{\chi} G_{ij}-\overline{\mathbb G^k_{\bar j}}G_{i\bar k},
\]
thus
\begin{equation}
\bar{\chi} G_{ij}
=\overline{\mathbb G^k_{\bar j}}G_{i\bar k}. \label{frak-D1'}
\end{equation}

%Set
%\[A:=\left(G_{ij}\right),\quad B:=\left(\mathbb G^i_{\bar j}\right).\]
%For K\"ahler--Finsler metric $F$ and fixed $(z,v)$, we can choose local normal coordinates such that $G_{i\bar j}(z,v)=\delta_{ij}$ by \cite{Ai91}. Under this coordinates, \eqref{frak-D1'} provides
%\begin{equation}
%\bar{\chi} A(z,v)
%=\overline{B}(z,v). \label{frak-D1}
%\end{equation}
Since $F$ is K\"ahler--Finsler, Lemma \ref{lem} gives
\begin{equation}
G_{ij}\mathbb G^i_{\bar l}=0 \quad\text{on}\quad \tilde M. \label{Ge-KF}
\end{equation}
%Since $A(z,v)$ is a symmetric matrix by definition, hence $B(z,v)$ is symmetric by \eqref{frak-D1}.
%Thus the transpose of \eqref{Ge-KF} gives
%\begin{equation}
%B(z,v)A(z,v)=\left(A(z,v)B(z,v)\right)^T=0.\label{van-symm}
%\end{equation}
Differentiate \eqref{Ge-KF} by $\bar{\chi}$ providing
\begin{equation}
0=\bar{\chi}\left(G_{ij}\mathbb G^i_{\bar l}\right)
=\bar{\chi}G_{ij}\mathbb G^i_{\bar l}
+G_{ij}\bar{\chi}\mathbb G^i_{\bar l}.\label{vanKF-1}
\end{equation}
Multiplication on the two sides of \eqref{vanKF-1} by $\mathbb G^j_{\bar m}$, take sum with $j$ and use \eqref{Ge-KF} providing
\begin{equation}
0=\bar{\chi}G_{ij}\mathbb G^i_{\bar l}\mathbb G^j_{\bar m}.\label{vanKF-2}
\end{equation}
Thus by \eqref{frak-D1'},
\begin{equation}
0=\overline{\mathbb G^k_{\bar j}}G_{i\bar k}\mathbb G^i_{\bar l}\mathbb G^j_{\bar m}.\label{vanKF-3}
\end{equation}
Multiplication on the two sides of \eqref{vanKF-3} by $\overline{\mathbb G^l_{\bar a}}$ and take sum with $l$ providing
\begin{equation}
0=\mathbb G^j_{\bar m}\overline{\mathbb G^k_{\bar j}}G_{i\bar k}\mathbb G^i_{\bar l}\overline{\mathbb G^l_{\bar a}}.\label{vanKF-4}
\end{equation}
Notice $\left(G_{i\bar k}\right)>0$, thus
\begin{equation}
0=\mathbb G^j_{\bar m}\overline{\mathbb G^k_{\bar j}}.\label{vanKF-5}
\end{equation}
By \eqref{frak-D1'},
\begin{equation}
0=\chi G_{\bar m\bar l}G^{\bar l j}\overline{\chi G_{\bar i\bar j}G^{\bar i k}}
=\chi G_{\bar m\bar l}G^{\bar l j}\bar\chi G_{ i j}G^{\bar k i}.\label{vanKF-6}
\end{equation}
Since $\left(G^{\bar k i}\right)>0$, we have
\begin{equation}
0=\chi G_{\bar m\bar l}G^{\bar l j}\bar\chi G_{ i j}.\label{vanKF-7}
\end{equation}
Since $\left(G^{\bar l j}\right)>0$, \eqref{vanKF-7} implies
\begin{equation}
0=\chi G_{\bar m\bar l}.\label{vanKF-7}
\end{equation}
Substitute \eqref{vanKF-7} into \eqref{frak-D1'} providing
\begin{equation}
0=\mathbb G^k_{\bar j}. \label{vanKF-8}
\end{equation}
Thus $\mathbb G^i$ is holomorphic in the fiber variables. Its $(2,0)$-homogeneity and smoothness on the slit fiber imply that it is quadratic:
\begin{equation}
\mathbb G^i=\mathbb G^i_{jk}(z)v^jv^k,\label{KF=KB1}
\end{equation}
where $\mathbb G^i_{jk}(z)
=\frac12\dot\partial_j\dot\partial_k\mathbb G^i(z,v)$.
By the (strongly) K\"ahler--Finsler condition, we have
$\Gamma^i_{j;k}=\Gamma^i_{k;j}$ and
\[\dot\partial_j\Gamma^i_{k;l}v^l
=\dot\partial_j\Gamma^i_{l;k}v^l
=\dot\partial_j\dot\partial_l\Gamma^i_{;k}v^l
=\dot\partial_l\Gamma^i_{j;k}v^l
=0,\]
where we used the Euler formula for the $(0,0)$-homogeneous map $\Gamma^i_{j;k}$
in the last equality. Thus
\begin{eqnarray*}
2\mathbb G^i_{jk}(z)
&=&\dot\partial_j\left(\Gamma^i_{k;l}(z,v)v^l
+\Gamma^i_{;k}(z,v)\right)\\
&=&\dot\partial_j\Gamma^i_{k;l}(z,v)v^l+\Gamma^i_{k;j}(z,v)
+\Gamma^i_{j;k}(z,v)\\
&=&2\Gamma^i_{j;k}(z,v).
\end{eqnarray*}
Thus $F$ is K\"ahler--Berwald.
\end{proof}

\section{Locally symmetric complex Finsler spaces}\label{sec-loc-symm}

A non-trivial map $\sigma: M \to M$ is called an involution at $p \in M$ if $\sigma(p) = p$ and $\sigma^2 = \operatorname{Id}$. And $\sigma$ is called a locally geodesic symmetry at $p$ if there exists a neighborhood $U \subset M $ of $p$, such that  $\sigma(p) = p$ and $\sigma(\gamma(t)) = \gamma(-t)$ for every real geodesic $\gamma \subset U$ passing through $p$ with $\gamma(0) = p$. We recall the following definitions.
\begin{definition}[ \cite{Deng}]\label{def2}
Let  $(M,F)$ be a complex Finsler manifold. Then $(M, F)$ (resp. $F$) is called a
\begin{itemize}
  \item symmetric complex Finsler space (resp. metric) if for every $p \in M$, there exists an involution of $M$ at $p$, which is a holomorphic isometry with $p$ as an isolated fixed point;
  \item locally symmetric complex Finsler space (resp. metric) if for every $p \in M$, there exists a neighborhood $U \subset M$ of $p$, and a geodesic symmetry of $M$ at $p$, which is a holomorphic isometry with $p$ as an isolated fixed point.
  \end{itemize}
\end{definition}
Let $\operatorname{Aut}(M)$ denote the group of holomorphic automorphisms of $M$. Homogeneous complex Finsler spaces are defined as follows.

\begin{definition}[\cite{Deng}]\label{def3}
Let $(M, F)$ be a connected complex Finsler manifold.
If for any points $p, q \in M$, there exists an isometry $f \in \operatorname{Aut}(M)$ such that $f(p) = q$. Then $(M, F)$ (resp. $F$) is called a homogeneous complex Finsler space (resp. metric).
\end{definition}

\begin{remark}
In  \cite{Deng}, S. Deng and Z. Hou proved that every symmetric (resp. homogeneous) complex Finsler space is underlyingly a symmetric (resp. homogeneous) Hermitian space. And analogous to the Hermitian case, one can prove that every locally symmetric complex Finsler space is homogeneous.
\end{remark}

The following lemma is similar to the Hermitian case, and proved with similar method.
\begin{lemma}\label{Lemma}
Let $(M, F)$ be a symmetric complex Finsler space and $p \in M$. If $\sigma$ is a holomorphic involutive isometry of $(M, F)$ such that $p$ is an isolated fixed point of $\sigma$, then $(d\sigma_p)\vert_p = -\operatorname{Id}.$
\end{lemma}
\begin{proof}
Let $A := (d\sigma_p)\vert_p$ and $M_p := T^{1,0}_p M$. Since $\sigma(p) = p$, we have $A M_p \subset M_p$. From $\sigma^2 = \operatorname{Id}$, it follows that $A^2 = \operatorname{Id}$. Consequently,
 \[A(v - Av) =(Av - A^2v)=-(v - Av),\quad A(v + Av)=(Av + A^2v)=(v+ Av).\]
Every $v \in M_p$ can be decomposed as $v = \frac{1}{2}(v - Av) + \frac{1}{2}(v + Av)$. Therefore, $M_p = V_- \oplus V_+$, where
\[V_-:= \{v:Av = -v\},\quad V_+:= \{v:Av = v\}.\]

Suppose there exists a non-zero $w \in V_+$. Consider the geodesic $\gamma$ through $p$ and with initial direction $w$. Since $\sigma$ is an isometry, $\sigma \circ \gamma$ is also a geodesic with the same initial direction $w$. By uniqueness, $\sigma$ fixes $\gamma$ pointwise, contradicting the assumption that $p$ is an isolated fixed point. Hence $V_+ = {0}$, so $M_p = V_-$ and $A = -\operatorname{Id}.$
\end{proof}

By the proof of Lemma \ref{Lemma}, for a symmetric complex Finsler space $(M, F)$ and $p \in M$, there exists a normal neighborhood $D_p$ of $p$, such that $\sigma_p$ acts as the geodesic symmetry on $D_p$. Thus we have the following corollary.
\begin{proposition}
 Symmetric complex Finsler spaces are locally symmetric  complex Finsler spaces.
\end{proposition}

From the proof of Lemma \ref{Lemma}, we obtain the following result.

\begin{lemma}\label{Lemma2}
Let $(M, F)$ be a locally symmetric complex Finsler space and $p \in M$. If $\sigma$ is a holomorphic geodesic symmetry of $(M, F)$ at $p$, then $(d\sigma_p)\vert_p = -\operatorname{Id}.$
\end{lemma}

Similar to the Hermitian case, we can obtain the following proposition for locally symmetric complex Finsler space.
\begin{proposition}
Let $M$ be a locally symmetric complex Finsler space, and $\mathcal G:=Iso(F)\bigcap Aut(M)$ be the group of biholomorphic maps of $M$ which are isometries of $F$. If $T$ is a s-form over $T^{1,0}M$ invariant under action of $\mathcal G$. Then
\begin{equation}
\sigma_p^\ast(T)=(-1)^sT,\label{tensor}
\end{equation}
for every isometric geodesic symmetry $\sigma_p$ with isolated fixed point $p\in M.$
\end{proposition}
\begin{proof}
Since $T$ is $\mathcal G$-invariant,$ a^\ast T= T,\; \forall \,a\in \mathcal G.$
Thus
\[
a^\ast(\sigma_p^\ast(T))=\sigma_p^\ast(T)\circ da=T\circ d(\sigma_p)\circ da,\,\forall \,a\in G.
\]
Hence $\sigma_p^\ast(T)$ is also a $\mathcal G$-invariant.
By Lemma \ref{Lemma2} we have $d(\sigma_p)\vert_p=-Id,$ thus
\[(\sigma_p^\ast(T))\vert_p=(-1)^sT\vert_p.\]
Since  $T$ and $(-1)^s\sigma_p^\ast(T)$ both are $\mathcal G$-invariance, and equal at $p$, thus they equal on $T^{1,0}M$.
\end{proof}

The K\"ahler form of complex Finsler metric $F$ is defined by
\[\omega:=\sqrt{-1}G_{i\bar j}(z,v)dz^i\wedge d\bar z^j.\]
Aikou proved that a complex Finsler metric is K\"ahlr--Finsler if and only if its K\"ahler form is closed under horizontal differential $d^{\mathcal H}$ \cite{Ai91}. Thus we obtain
\begin{proposition}\label{Kahler-F}
Strongly pseudoconvex locally symmetric complex Finsler metrics are K\"ahler--Finsler.
\end{proposition}
\begin{proof}
Since $\sigma_p$ is holomorphic, the actions of $d^{\mathcal H}$ and $\sigma_p^\ast$ are  interchangeable, we have
\begin{equation}
d^{\mathcal H}(\sigma_p^\ast(\omega))=\sigma_p^\ast(d^{\mathcal H}(\omega)).\label{d-sigma}
\end{equation}
Thus by \eqref{tensor}, the left and right hand side of \eqref{d-sigma} are $
(-1)^kd^{\mathcal H}(\omega)$ and $(-1)^{k+1}d^{\mathcal H}(\omega),$ respectively,
which implies $d^{\mathcal H}\omega=0.$
\end{proof}

By Theorem \ref{thm:KF=KB} and Proposition \ref{Kahler-F}, we have
\begin{theorem}\label{KF=symm}
Strongly pseudoconvex locally symmetric complex Finsler metrics are K\"ahler--Berwald.
\end{theorem}

In order to study the relationship between locally symmetric complex Finsler metrics and locally symmetric Hermitian metrics, we provide another proof of Theorem \ref{KF=symm} bolow. Indeed, we obtain the following proof more earlier than the above one, which is independent of Theorem \ref{Kob-KB} and results in \cite{Ge-wKF}.

We recall at first some known results and property about real Berwald and complex (resp. K\"ahler--)Berwald metrics,
in particular results of Aikou \cite{Ai95}, which can help connecting complex Berwald (resp. K\"ahler--Berwald) metrics to Hermitian (resp. K\"ahler) ones.

\begin{theorem}[ \cite{Ai95}, Theorem 1.1]\label{Ai1.1}
 Let $E$ be a smooth complex vector bundle over $M$ with a complex Finsler metric $F(x, \xi)$ . We suppose that $E$ admits a complex connection $\nabla$ on $E$ which preserves the norm invariant under the parallel displacement. Then there exists a Hermitian metric $h$ on $E$ such that $\nabla$ is a metrical connection of $(E, h)$ .
\end{theorem}
The Hermitian metric $h$ is induced from $H_p^0$-invariant Hermitian inner product $\langle\cdot,\cdot\rangle$ var the parallel displacement with respect to $\nabla$:
\[
h(\xi,\psi):=
\left \langle P_{c}^{-1} \xi, P_{c}^{-1} \psi\right\rangle_{p},\quad \forall \,\xi,\psi\in E_q,
\]
where $H_p^0$ is unitary component of the holonomy group of $\nabla$ with reference point $p \in M$, $c(t) (0 \leq t\leq 1)$ is a $C^{\infty}$ curve such that $c(0)=p$ and $c(1)=q$, and $P_{c}: E_{p} \to E_{q}$ is the parallel displacement with respect to $\nabla$ along $c(t)$.
\begin{remark}\label{equi-unifm}
Such a $H_p^0$-invariant inner product can be obtained from the following method from any inner product $(\cdot,\cdot)$:
\[
\langle u,v\rangle:=\int_{\mathcal G} (g(u),g(v)) dg,
\]
where $dg$ is the bi-invariant Haar measure on $\mathcal G,$ and $\mathcal G$ is the group invariant under the Finsler norm $F(p;\cdot):$
\[\mathcal G:=\{g\in Gl(n;\mathbb C): F(p,v)=F(p,g(v)),\quad\forall \, v\in T^{1,0}_pM\}.\]
$\langle \cdot,\cdot\rangle$ is invariant under  $\mathcal G$, hence invariant under the subgroup $H_p^0< \mathcal G$.

Since $F(\xi)= F(P_{c}^{-1} \xi),\quad \forall \,\xi\in E_q$ by condition, thus $F$ and $h$ are both invariant under the parallel displacement $P_c$.
Since any two continuous norms are equivalent on a vector space,
thus $F$ and $h$ are equivalent uniformly.
\end{remark}
Aikou's  Theorem 4.1, Propositiion 4.3 and Proposition 4.4 in \cite{Ai95} imply the following results.
\begin{theorem}[\cite{Ai95}]\label{Ai-thm}
 Let $(M, F)$ be a complex Berwald manifold. Then there exists a Hermitian metric $h_{M}$ on $M$ such that the Chern--Finsler nonlinear connection of $F$ coincides with the Chern connection of $h_{M}$.
And $F$ is a K\"ahler--Finsler metric if and only if $h_{M}$ is a K\"ahler metric.
\end{theorem}
\begin{remark}
Here the pull-back bundle is $\pi^\ast(T^{1,0}M)$ over the split bundle $\tilde M$. And Aikou considered the Chern--Finsler connection as
 $\nabla:\Gamma(\pi^\ast(T^{1,0}M))\rightarrow \Gamma(T^\ast_C\tilde M\otimes \pi^\ast(T^{1,0}M)).$
However, the Chern--Finsler connection associated to the vertical or horizontal sub-bundle $\mathcal V$ and $\pi^\ast(T^{1,0}M)$ coincide excepted for their vector bases.
\end{remark}

Before proving the main theorem of this section, let us recall some notations and result of S. Deng.
A holomorphic (resp. diffeomorphism) map $f:(M,\nabla)\rightarrow (N,\hat\nabla)$ is called an affine transformation if
\[
\hat\nabla_{df (X)}df (Y)=df(\nabla_XY),
\]
for any $X\in T^{1,0}M$  (resp. $X\in TM$) and smooth section $Y$ of $T^{1,0}M$ (resp. $TM$).
S. Deng find that the relationship between affine translations and isometries for Riemann geometry also holds in the real Berwald setting. % (cf. Lemmea 1.4 of Section 1, Chapter IV in \cite{Heg}) into real Finsler circumstance.
\begin{theorem}[ \cite{Deng2}, Theorem 5.1]\label{Deng5.1}
 Let $(M,F)$ be a connected real Berwald space and $\psi$ an isometry of $(M,F)$ onto itself. Then $\psi$ is an affine transformation with respect to the connection of $F$.
\end{theorem}

With similar methods used in the real Berwald case, Theorem \ref{Deng5.1} also holds in the complex Berwald setting.
Before the proof, we introduce the following induced Hermitian metric $g_w$ on $M$ using fixed holomorphic vector field $w\in \Gamma\left(T^{(1,0)}M\right)$:
\[g_w(u,v):=G_{i\bar j}(z,w(z))u^i\bar v^j, \, u,v\in T^{(1,0)}_zM.\]
Clearly, with complex parameters $s,t\in \mathbb C$, it holds that
\begin{equation}
g_w(u,v)=\frac{\partial^2}{\partial  s\partial \bar t}\vert_{s=t=0}G(z,w+su+tv).\label{gw}
\end{equation}

\begin{lemma}\label{KB-Deng5.1}
 Let $(M,F)$ be a connected complex Berwald space and $f$ a locally biholomorphic isometry. Then $f$ is an affine transformation with respect to the Chern-Finsler nonlinear connection of $F$:
\[\nabla_{df(u)}df(w)
=df\left(\nabla_{u}w\right),\,\forall \,w\in\Gamma(T^{(1,0)}M),u\in T_p^{(1,0)}M.\]
\end{lemma}
\begin{proof}
Similar to real Finsler case (see page 260 of \cite{Bao-Chern-Shen}), take locally holomorphic curve $\gamma:\epsilon \mathbb D\rightarrow M$ with initial value
$\gamma(0)=p,\,\frac{\partial\gamma}{\partial \zeta}(0)=u,$ then
\begin{eqnarray}
u[g_w(w,v)]&=& \frac{\partial }{\partial \zeta}\vert_{\zeta=0}\left(G_{i\bar j}(\gamma,w\circ\gamma)\right)w^i(p)\bar v^j(p)+G_{i\bar j}(p,w(p))u(w^i\bar v^j).\label{d0}
\end{eqnarray}
Notice $G_{i\bar jl}(z,w)w^i=0,$ we have
\begin{eqnarray*}
&&\frac{\partial }{\partial \zeta}\vert_{\zeta=0}\left(G_{i\bar j}(\gamma,w\circ\gamma)\right)w^i(p)\\
&=&\left(G_{i\bar j;k}(p,w(p))
+G_{i\bar jl}(p,w(p))\frac{\partial w^l}{\partial \gamma^k}(0)\right)u^kw^i(p)\\
&=&\left(G_{i\bar j;k}-\Gamma^l_{;k}G_{i\bar jl}\right)
(p,w(p))u^kw^i(p)\\
&=&\delta_kG_{i\bar j}(p,w(p))u^kw^i(p)\\
&=&G_{l\bar j}(p,w(p))\Gamma^l_{i;k}(p,w(p))u^kw^i(p)
\end{eqnarray*}
Substitute this into \eqref{d0}, we have
\begin{eqnarray}
u[g_w(w,v)]&=&g_{w(p)}\left(\nabla_uw, v(p)\right)
+g_{w(p)}\left(w(p),\bar u(v^j)\frac{\partial}{\partial z^j}\right).\label{d1}
\end{eqnarray}
By equation \eqref{gw}, for any biholomorphic isometric map $f$ of $(M,F)$, we have
\begin{eqnarray}
g_{df(w)}(df(w),df(v))=g_w(w,v).\label{d2}
\end{eqnarray}
By \eqref{d1} we have
\begin{eqnarray}
&&g_{df(w(p))}\left(\nabla_{df(u)}df(w),df(v(p))\right)\\
&=&df(u)g_{df(w)}(df(w),df(v))
-g_{df(w(p))}\left(df(w(p)),\overline{df(u)}
(df^i(v))\frac{\partial}{\partial z^i}\right).\label{d3}
\end{eqnarray}
Substitute \eqref{d2} into the right hand of \eqref{d3},
\begin{eqnarray}
g_{df(w(p))}\left(\nabla_{df(u)}d\psi(w),df(v(p))\right)
&=&ug_{w}(w,v)-g_{w(p)}\left(w(p),\overline{u}(v^i)\frac{\partial}{\partial z^i}\right)\nonumber\\
&=&g_{w(p)}\left(\nabla_uw,v(p)\right),\label{d4}
\end{eqnarray}
where in the last equation we used \eqref{d1}.
Substitute \eqref{d2} into \eqref{d4}, we have
\begin{eqnarray}
g_{df(w(p))}\left(\nabla_{df(u)}df(w),df(v(p))\right)
&=&g_{df(w(p))}\left(df\left(\nabla_{u}w\right),
df(v(p))\right),\label{d5}
\end{eqnarray}

Since $g_{df(w(p))}$ is a Hermitian inner product, by the arbitiary of $v$, \eqref{d5} implies
\[\nabla_{df(u)}df(w)=df\left(\nabla_{u}w\right).\]
\end{proof}

Now we can obtain the main theorem of this section.

\begin{theorem}\label{thm1}
Every strongly pseudoconvex locally symmetric complex Finsler metric  $F$ is a K\"ahler--Berwald, complex Finsler--Einstein metric. And the corresponding Hermitian metrics whose Chern connection coincide with the Chern--Finsler nonlinear connection of $F$ are locally symmetric.
\end{theorem}
\begin{proof}
 For any $p\in M$, by definition, there exists a biholomorphic mapping $f$ of $M$, which is an isometric geodesic symmetry in a neighborhood $U$ of $p$ with $p$ as an isolated fixed point. In the coordinate system $(z,v)$ of $T^{1,0}U$, we have $F(z,v)=F(f(z);df(v)).$ Denote the Jacobi matrix $(\frac{\partial f^i}{\partial z^j})$ by $(J^i_j)$, thus
 \begin{eqnarray*}
G_{\bar b}(z,v)=G_{\bar\beta}(f(z),df(v))\bar J^\beta_b,\label{fz1}
 \end{eqnarray*}
 \begin{eqnarray}
 G_{\bar b;c}(z,v)=G_{\bar\beta;\xi}(f(z),df(v))\bar J^\beta_b J^\xi_c+G_{\bar\beta\xi}(f(z),df(v))\bar J^\beta_b\frac{\partial^2f^\xi}{\partial z^a\partial z^c}(z)v^a.\label{fz2}
\end{eqnarray}
 Since $p$ is the isolated fixed point of $f$, by Lemma \ref{Lemma2} we have
 \begin{equation}
 f(p)=p,df_p=-Id.\label{-Id}
 \end{equation}
Hence restricted at $p$, \eqref{fz2} can be reduced to
 \begin{eqnarray}
 G_{\bar b;c}(p,v)&=&G_{\bar b;c}(p,-v)-G_{\bar b\xi}(p,-v)\frac{\partial^2f^\xi}{\partial z^a\partial  z^c}(p)v ^a\nonumber\\
 &=&-G_{\bar b;c}(p,v)-G_{\bar b\xi}(p,v)\frac{\partial^2f^\xi}{\partial z^a\partial  z^c}(p)v^a,\nonumber\\
\end{eqnarray}
where in the last equality we used the homogeneous of $G_{\bar b;c}$ and $G_{\bar b\xi}.$ Thus
 \begin{eqnarray}
G_{\bar b;c}(p,v)&=&-\frac12G_{\bar b\xi}(p,v)\frac{\partial^2f^\xi}{\partial z^a\partial  z^c}(p)v^a.
 \label{fz4}
\end{eqnarray}
 By \eqref{fz4} we have
 \begin{eqnarray}
 \Gamma^a_{;c}(p,v)
 &=&G^{\bar ba}(p,v)G_{\bar b;c}(p,v)\nonumber\\
 &=&-\frac12G^{\bar ba}(p,v)G_{\bar b\xi}(p,v)\frac{\partial^2f^\xi}{\partial z^d\partial \bar z^c}(p)v^d\nonumber\\
 &=&-\frac12\frac{\partial^2f^a}{\partial z^d\partial \bar z^c}(p)v^d.
 \end{eqnarray}
 Thus
 \begin{eqnarray}
  \Gamma^a_{b;c}(p,v)
 =\frac{\partial}{\partial v^b}\Gamma^a_{;c}(p,v)
 =-\frac12\frac{\partial^2f^a}{\partial z^b\partial  z^c}(p).\label{fz5}
 \end{eqnarray}
By \eqref{fz4}, the horizontal Chern--Finsler connection coefficients of $F$
are independent of $v$, and symmetric with respect to lower index, thus $F$ is K\"ahler--Berwald.

By Theorem \ref{Ai-thm}, there exists a K\"ahler metric $h$, whose Chern connection coincides with the Chern--Finsler nonlinear connection of $F$, and thay share same geodesics as parameterized curves.
Thus $f$ is an geodesic symmetry of $(M,h)$ with $p$ as an isolated fixed point.
By Lemma \ref{KB-Deng5.1}, $f$ is an affine locally biholomorphic transformation of $(M,F)$ with respect its Chern--Finsler nonlinear connection, thus an  affine locally biholomorphic transformation of $(M,h)$ with respect to its Chern connection.
Since
\[h(p,v)=h(p,-v)=h(f(p),df\vert_p(v)),\]
$f$ is an isometry of $(M,h)$ when restricted at $T^{1,0}_pM,$ thus an isometry of $(M,h)$ by Lemmea 1.4 of Section 1, Chapter IV in \cite{Heg}. Thus $h$ is a locally symmetric K\"ahler metric.
Since a locally symmetric Hermitian metric is K\"ahler--Einstein, by Theorem 3.5 of \cite{SZ18}, the locally symmetric complex Finsler metric $F$ is a complex Einstein--Finsler metric in the sense of Aikou \cite{Ai97}.
\end{proof}

The following characterization of locally symmetric complex Finsler spaces is a generalization of the Hermitian case.
\begin{theorem}\label{pal-eq}
A complete strongly pseudoconvex K\"ahler--Berwald metric $F$ is locally symmetric if only if
\[\nabla_H R^{\mathcal H} =0,\quad\forall \,H\in \mathcal H,\]
where $R^{\mathcal H}:=R_{j \bar s;k\bar l}dz^k\wedge d\overline{z^l}\otimes dz^j \otimes d\overline{z^s},$ and $R_{j \bar s;k\bar l}$ is the $h\bar h$-part Chern--Finsler curvature defined by \eqref{HH-curv'}.
\end{theorem}
\begin{proof}
If $F$ is a strongly pseudoconvex locally symmetric complex Finsler mettric, by Theorem \ref{thm1}, its nonlinear connection $\nabla$ coincides with the Chern connection of a locally symmetric Hermitian metric $h$. By \eqref{HH-curv'},
\begin{equation}
R_{j \bar s;k\bar l}=G_{i\bar s}R_{j;k\bar l}^i,\,
R_{j;k\bar l}^i=-\delta_{\bar l}\Gamma^i_{j;k}
-\Gamma^i_{js}\delta_{\bar l}\Gamma^s_{;k},\label{H-curv}
\end{equation}
and the complex Berwald condition makes
\[R_{j;k\bar l}^i=-\partial_{\bar l}\Gamma^i_{j;k},\]
which is also the Chern curvature of $h$.  Since $h$ is K\"ahler, $\nabla$ is also the Levi--Civita connection, by Cartan's characterization of locally symmetric space, the Chern curvature $R^h$ of $h$ is parallel with respect to $\nabla$, i.e., for $\forall \,v\in \tilde M,$
\[
\begin{aligned}
0&=\nabla_{\partial_m} R^h\\
&=\nabla_{\partial_m} \left(R^i_{j;k\bar l}dz^k\wedge d\overline{z^l}\otimes dz^j \otimes h_{i\bar s} d\overline{z^s}\right)\\
&=\nabla_{\partial_m} \left(R^i_{j;k\bar l}dz^k\wedge d\overline{z^l}\otimes dz^j \right)\otimes h_{i\bar s} d\overline{z^s}+R^i_{j;k\bar l}dz^k\wedge d\overline{z^l}\otimes dz^j \otimes\nabla_{\partial_m} \left( h_{i\bar s} d\overline{z^s}\right)
 .\end{aligned}
 \]
That is
\begin{equation}
\nabla_{\partial_m} \left(R^i_{j;k\bar l}dz^k\wedge d\overline{z^l}\otimes dz^j\right)
=-\Gamma_{\tau;m}^{i}R^\tau_{j;k\bar l}dz^k\wedge d\overline{z^l}\otimes dz^j.\label{parallel-curvature}
\end{equation}
 Notice
 \begin{eqnarray}
 \nabla_{\delta_m} R^{\mathcal H}
& =&\nabla_{\delta_m} \left(R^\alpha_{\beta;\gamma\bar \delta}dz^\gamma\wedge \overline{dz^\delta}\otimes dz^\beta\right) \otimes G_{\alpha\bar\sigma } \overline{dz^\sigma}\nonumber \\ &&+R^\alpha_{\beta;\gamma\bar \delta}dz^\gamma\wedge \overline{dz^\delta}\otimes dz^\beta \otimes \nabla_{\delta_m} \left(G_{\alpha\bar\sigma } \overline{dz^\sigma}\right)\nonumber\\
& =&\nabla_{\partial_m} \left(R^\alpha_{\beta;\gamma\bar \delta}dz^\gamma\wedge \overline{dz^\delta}\otimes dz^\beta\right) \otimes G_{\alpha\bar\sigma } \overline{dz^\sigma}\nonumber \\ &&+R^\alpha_{\beta;\gamma\bar \delta}dz^\gamma\wedge \overline{dz^\delta}\otimes dz^\beta \otimes G_{\tau\bar\sigma } \Gamma^\tau_{\alpha;m}\overline{dz^\sigma}, \label{parallel-curvature1}
 \end{eqnarray}
where in the last equality, we use \[\delta_mR^\alpha_{\beta;\gamma\bar \delta}(z)
=\partial_mR^\alpha_{\beta;\gamma\bar \delta}(z)
\]
and the horizontal-metrical property of the Chern--Finsler connection:
\[\delta_m(G_{\alpha\bar\beta})
=G_{\tau\bar\beta}\Gamma^\tau_{\alpha;m}
+G_{\alpha\bar\tau}\overline{\Gamma^\tau_{\beta;m}}.
\]
Substitute \eqref{parallel-curvature} to \eqref{parallel-curvature1} provides
$\nabla_{\delta_m} R^{\mathcal H}=0.$

On the other hand, if
\[\nabla_H R^{\mathcal H}=0, \quad \forall \,H\in \mathcal H,\]
then
 \begin{eqnarray*}
0= \nabla_{\delta_m} R^{\mathcal H}
& =&\nabla_{\partial_m} \left(R^\alpha_{\beta;\gamma\bar \delta}dz^\gamma\wedge d\overline{z^\delta}\otimes dz^\beta\right) \otimes G_{\alpha\bar\sigma } d\overline{z^\sigma} \\ &&+R^\alpha_{\beta;\gamma\bar \delta}dz^\gamma\wedge d\overline{z^\delta}\otimes dz^\beta \otimes G_{\tau\bar\sigma } \Gamma^\tau_{\alpha;m}d\overline{z^\sigma},
 \end{eqnarray*}
 thus \eqref{parallel-curvature} holds, hence $\nabla_{\partial_m} R^h=0$. Since $h$ is K\"ahler, $R^h$ is real tensor, we also have
 $\nabla_{\partial_{\bar m}} R^h=0$, thus $\nabla R^h=0$.
Since $\nabla$ is also the Levi--Civita connection,
$h$ is locally symmetric.
By the proof of Theorem \ref{thm1}, $F$ is also a locally symmetric complex Finsler metric.
\end{proof}

In \cite{Deng-Hou05},
S. Deng and Z. Hou obtained the flag curvature of strongly convex real locally symmetric real Berwald spaces, and obtained the real flag curvature properties for real locally symmetric Finsler spaces.
Even though a strongly pseudoconvex complex Finsler metric need not to be a  strongly convex real Finsler metric,
by Theorem \ref{thm1}, strongly pseudoconvex locally symmetric complex Finsler metrics are K\"ahler--Berwald metrics, thus we can also associated their curvature with corresponding Hermitian symmetric metrics, thus we can obtain similar curvature formulas and their signs as in \cite{Deng-Hou05}.
The proof is similar, we represent it here for completeness. Here we denote $M=Q/H$, where $Q$ is the unitary connect component of the holomorphic isometric group of $(M,F)$, and $H$ the isotropy subgroup of $Q$ at a fixed point $z\in M,$ $\mathfrak q$ and $\mathfrak h$ be the Lie algebras of $Q$ and $H$, respectively, and $\mathfrak q=\mathfrak h\oplus \mathfrak m$ be the canonical decomposition of the Cartan symmetric pair $(\mathfrak q,\mathfrak h)$.

\begin{proposition}\label{3curvature}
 Let $(M=Q/H, F)$ be a locally symmetric complex Finsler spaces of compact type. Then the curvature tensor of $F$ is given by
\[\Omega(\bar u,v)w=[[v,\bar u],w], \forall \,u,v,w\in \mathfrak m.\]
The holomorphic sectional and bisectional curvatures of $F$ at $o$ are given respectively by
	\begin{equation}
		K(o,v)=\frac{2}{G^2(o,v)}\langle \chi_v\left([[v,\bar v]v]\right),\chi_v(v)\rangle_{(o,v)}
,\quad v\in \mathfrak m;\label{4hsc}
	\end{equation}
\begin{equation}
		B(o,v,w)=\frac{2}{G(o,v)G(o,w)}\left\langle \chi_v\left([[v,\bar w]w]\right),\chi_v(v)\right\rangle_{(o,v)},
\quad v,w\in \mathfrak m,\label{4hbsc}
	\end{equation}
where $\chi_v(w)=\sum_{i=1}^nw^i\delta_i\vert_{(o;v)}$ is the horizontal lifting of $w$ along $v$.
While these curvatures differ by a negative sign for noncompact type symmetric spaces, they vanish for Euclidean type symmetric spaces.
Thus $K$ is bounded between negative constants (resp. between positive constants, or identically zero) if $M$ is of non-compact type (resp. compact type, or Euclidean type),
 $B$ is semi-negative (resp. semi-positive, or zero) if $M$ is of non-compact type (resp. compact type, or Euclidean type). Furthermore, $B$ has a negative lower bound (resp. a positive upper bound) if $M$ is of non-compact type (resp. compact type).
\end{proposition}
\begin{proof}
By Theorem \ref{thm1}, $(Q/H, F)$ is a globally symmetric K\"ahler--Berwald space, and the horizontal part of the (linear) Chern--Finsler connection of $(Q/H, F)$ $\nabla$ is just the Chern connection of the globally symmetric Hermitian space $(Q/H,h)$, therefore, for compact-type case (cf. \cite{Borel60}),
\[\Omega(\bar u, v)w=[[v,\bar u],w],\quad u,v,w\in \mathfrak m.\]
And \eqref{4hsc}-\eqref{4hbsc} follows from their definitions.

When $M$ has complex dimension $1$, every complex Finsler metric on $M$ is Hermitian, thus the results follows. Otherwise, we always have $dim_{\mathbb R}\mathfrak h\geq2,$ thus $\mathfrak h$ contains complex subspace. By the invariance of $F$ under $Ad(H)$, the adjoint representation of $G$ when restricted on $H$ (i.e., the  isotropy representation), for $t\in \epsilon\mathbb D\subset \mathfrak h$ we have
 \[G(Ad(exp(tw))(y+ru+sv))=G(y+ru+sv),\quad w\in \mathfrak h,\, r,s\in \mathbb C,\, u,v,0\neq y\in \mathfrak m.\]
Since \[ Ad(exp(tw))u=u+t[w,u]+O(\vert t\vert^2),\quad u\in \mathfrak m, \]
thus
\begin{eqnarray*}
\langle u,v\rangle_y
&=&\frac{\partial^2G(y+ru+sv+t[w,y+ru+sv]+O(\vert t\vert^2))}{\partial r\partial \bar s}\vert_{r=s=0}\\
&=&\langle u+t[w,u],v+t[w,v]\rangle_{y+t[w,y]}.
\end{eqnarray*}
Taking derivative with respect to $t$ leads to
\begin{equation}
0=\langle [w,u],v\rangle_{y}+\langle u,[w,v]\rangle_{y}
+G_{\alpha\bar\beta\gamma}(y)[w,y]^\gamma u^\alpha \bar v^\beta.\label{cCartan}
\end{equation}
Take $ u=y$ in \ref{cCartan}, we have
\[0=\langle [w,y],v\rangle_{y}+\langle y,[w,v]\rangle_{y}.\]
Which implies
\[\frac{\partial}{\partial t}\langle Ad(exp(tw))y, Ad(exp(tw))v\rangle_y=0,\]
notice $H$ is connected and compact, the exponent map is surjective, thus
\[\langle Ad(h)y, Ad(h)v\rangle_y=\langle y, v\rangle_y,\quad \forall \,h\in H.\]
For each fixed $y\in \mathfrak m$, we can define a $H$-invariant inner product on $\mathfrak m$ by
\[\langle u,v\rangle_p:=\int_{H} \langle Ad(h)u, Ad(h)v\rangle_y dh,\quad \forall \,u,v\in \mathfrak m,\]
where $dh$ is the normalized standard bi-invariant Haar measure of $H$, thus
\[\langle y,v\rangle_p=\langle y,v\rangle_y.\]
And we can extend $ \langle \cdot,\cdot \rangle_p$ to a $ Q$-invariant Hermitian metric $g$ on $Q/K,$ such that $(Q/K,g)$ is a Hermitian symmetric space. Since $g$ and $F$ have same horizontal Chern--Finsler connection, and by  \eqref{4hsc}-\eqref{4hbsc}, the resulds hold for the fixed vector $y$. By the arbitrary of $y\in \mathfrak m$, we have done.
\end{proof}

\section{de Rham decomposition theorem}\label{sec-de-Rham}

 In this section, we give the de Rham decomposition theorem  for strongly pseudoconvex K\"ahler--Berwald manifolds. Then we verify the generalized Frankel conjecture in this setting.

The Descarts product for K\"ahler--Berwald manifold is defined similar to real Berwald case, and it is a generalization of the Descarts product for K\"ahler manifold \cite{Zh2}. We only introduce for product of two manifolds, general cases are similar.

For a product manifold $M^{(1)}\times M^{(2)}$ of Hermitian manifolds, with projection map $\pi^{(i)}:M^{(1)}\times M^{(2)}\rightarrow M^{(i)},i=1,2.$  One has the direct sum
\[
T^{1,0}_p\left(M^{(1)}\times M^{(2)}\right)= T^{(1)}_p\oplus T^{(2)}_p, \,\forall \,p\in M^{(1)}\times M^{(2)},
\]
such that $\pi_\ast^{(1)}T^{(2)}_p=0$ and  $\pi_\ast^{(2)}T^{(1)}_p=0.$
For $v^{(i)}\in T^{1,0}M^{(j)},$ let $\pi^{\ast0}(v^{(i)})\in T^{1,0}_pM^{(1)}\times M^{(2)}$ be the unique vector such that
\[\pi^{(i)}_\ast\pi^{\ast0}(v^{(i)})=v^{(i)},\,
\pi^{(j)}_\ast\pi^{\ast0}(v^{(i)})=0,\,i\neq j.
\]
The complex linear connection $\nabla^{(1)}\times \nabla^{(2)}$ on $M^{(1)}\times M^{(2)}$ denotes the Descarts product of  complex linear connections $\nabla^{(i)}$ on $M^{(i)},\,i=1,2,$ i.e.,
\[\nabla^{(1)}\times \nabla^{(2)}_{\pi^{\ast0}(v^{(i)})}\pi^{\ast0}(w^{(i)})
=\pi^{\ast0}\left(\nabla^{(i)}_{v^{(i)}}w^{(i)}\right),\, v^{(i)},w^{(i)}\in \Gamma(T^{1,0}M^{(i)}),\]
and
\[\nabla^{(1)}\times \nabla^{(2)}_{\pi^{\ast0}(v^{(j)})}\pi^{\ast0}(v^{(i)})=0,\, v^{(k)}\in \Gamma(T^{1,0}M^{(k)}),k=1,2,i\neq j.\]

\begin{definition}[\cite{Zh2}]
A  complex Berwald space $(M,F,\nabla)$ is side to be the Descartes product of two complex Berwald spaces $(M^{(i)},F^{(i)},\nabla^{(i)}),i=1,2,$ if
\begin{enumerate}
  \item $(M,\nabla) = (M^{(1)}\times M^{(2)},\nabla^{(1)}\times\nabla^{(2)}),$
  \item $F(\pi^{\ast0}(v^{(i)}))=F^{(i)}(v^{(i)}),\,\forall \, v^{(i)}\in T^{1,0}(M^{(i)})$,
\end{enumerate}
where $\nabla$ and $\nabla^{(i)}$ are the Chern--Finsler nonlinear connection of $F$ and $F^{(i)},i=1,2$, respectively.
\end{definition}
\begin{remark}
The conditions does not imply that $G$ ia a product metric. For example, on a Hermitian manifold $(M,h)=\left(M^{(1)}\times M^{(2)},h^{(1)}\times h^{(2)}\right), $ the complex Finsler metric
\[F^2(z,v)=\sum_{i=1}^2h^{(i)}(z^{(i)},v^{(i)})
+t\sqrt[k]{\sum_{i=1}^2\left(h^{(i)}(z^{(i)},v^{(i)})\right)^k}, (z^{(i)},v^{(i)})\in T^{1,0}M^{(i)}\]
is a strongly convex complex Berwald metric, and $(M,F,\nabla)$ is Descarts product of $(M^{(1)},h^{(1)},\nabla^{(1)})$ and $(M^{(2)},h^{(2)},\nabla^{(2)})$ for any fixed $t>0$ and $k\in\mathbb N_+$, and $F$ is a K\"ahler--Berwald metric if $h^{(j)},\, j=1,2,$ are K\"ahler metrics (cf. \cite{Zh2}).
\end{remark}

In our proof, the following theorem of Berger and Simons about the holonomy group of Riemann (especially, K\"ahler) metric is needed \cite{Ber55,Sim62}.

\begin{theorem}[\cite{Ber55,Sim62}]\label{Berger}
Let $(M,g)$ be an irreducible Riemannian manifold. If $Hol_p^o$ is not transitive on the unit sphere in $(T_pM,g(p,\cdot))$, then $M$ is a symmetric space of rank $\geq2$.
\end{theorem}
\begin{remark}
Here $Hol_p^o$, the unitary component of the holonomy group of $(M,g)$, is with respect to an affine connection with zero-torsion $\hat\nabla$. And the curvature  $\hat R$ of $\hat\nabla$ is parallel with respect to the connecction. Thus when  $\hat\nabla$ is the Riemann connection of a Riemann metric $g$, $(M,g)$ is a Riemann symmetric space.
\end{remark}

Before proving of the de Rham decomposition theorem, we generalize a result for an affine transformation to be isometry in the complex Berwald setting.
Which was similar to the real Berwald case (cf. \cite[Theorem 5.2]{Deng2}).

\begin{theorem}[\cite{Deng2}]\label{aff}
 Let $(M_i,F_i)$, $i= 1,2,$ be two real Berwald spaces, and $\psi$ an affine diffeomorphism from $(M_1,F_1,\nabla_1)$ onto $(M_2,F_2,\nabla_2)$. If there exists $p\in M_1$ such that $d\psi$ is a linear isometry from $T_p(M_1)$ onto $T_{\psi(p)}(M_2)$, then $\psi$ is an isometry.
\end{theorem}
The following lemma can be proved with similar methods of Theorem \ref{aff}.
\begin{lemma}\label{aff-iso}
 Let $(M_i,F_i)$, $i= 1,2,$ be two complex Berwald spaces, and $f$ an affine locally biholomorphism from $(M_1,F_1)$ onto $(M_2,F_2)$ with respect to the Chern--Finsler nonlinear connections of $F_1$ and $F_2$. If there exists $p\in M_1$ such that $d\psi$ is a linear isometry from $T_p^{1,0}(M_1)$ onto $T_{f(p)}^{1,0}(M_2)$, then $f$ is a locally isometry.
\end{lemma}

\begin{proof}
Let $q \in M_{1}$, we need to prove $F_{1}(y)=F_{2}(df(y))$, $\forall \,y \in T_{q}^{1,0}(M_{1})$.

\textbf{Step 1.}
Join $q$ to $p$ by a piecewise smooth curve $\gamma:(-2,2)\rightarrow M,\, \gamma(0)=q,\gamma(1)=p.$
For every $u,v\in T^{1,0}_qM$, we still denote by $y,u,v$ the parallel vector fields along $\gamma$. Define
\[g(t):=g_{y(t)}(u(t),v(t))=G_{a\bar b}(\gamma(t),y(t))u^a(t)\bar v^b(t).\]
Below we prove that $g(t)$ is constant, and we use the notation $\dot \varphi:=\frac{d}{dt}\varphi$ in the following calculation.
\begin{eqnarray}
\dot g
&=&\Big[ G_{a\bar b;c}(\gamma,y)\dot\gamma^c
+G_{a\bar b;\bar c}(\gamma,y)\dot{\bar\gamma}^c
+G_{a\bar bc}(\gamma,y)\dot y^c+G_{a\bar b\bar c}(\gamma,y)\dot{\bar y}^c\Big]u^a\bar v^b\nonumber\\
&&+G_{a\bar b}(\gamma,y)\Big[\dot u^a\bar v^b+u^a\dot{\bar v}^b\Big].\label{affiso0}
\end{eqnarray}
By complex Berwald condition, $\Gamma^a_{b;c}(\gamma,\cdot)=
\Gamma^a_{b;c}(\gamma)$. And by the equation of parallel vector field for complex Berwald metrics, \eqref{geod},
\begin{equation}
\dot y^a=-\Gamma^a_{b;c}(\gamma)y^b\dot\gamma^c,\,
\dot u^a=-\Gamma^a_{b;c}(\gamma)u^b\dot\gamma^c,\,
\dot v^a=-\Gamma^a_{b;c}(\gamma)v^b\dot\gamma^c.
\label{paral-v}
\end{equation}
Substitute \eqref{paral-v} into \eqref{affiso0}, and then by $\Gamma^c_{d;e}y^d=\Gamma^c_{;e}(y)$, we have
\begin{eqnarray}
\dot g
&=&\Big( G_{a\bar b;c}(\gamma,y)\dot\gamma^c
+G_{a\bar b;\bar c}(\gamma,y)\dot{\bar\gamma}^c
-G_{a\bar bc}(\gamma,y)\Gamma^c_{d;e}(\gamma)y^d\dot\gamma^e
+G_{a\bar b\bar c}(\gamma,y)\nonumber\\
&&\overline{\Gamma^c_{d;e}(\gamma)y^d\dot{\gamma}^e}\Big)
u^a\bar v^b
-G_{a\bar b}(\gamma,y)\Big(\Gamma^a_{d;e}(\gamma)u^d\dot\gamma^e
\bar v^b+u^a\overline{\Gamma^b_{d;e}(\gamma)v^d\dot{\gamma}^e}\Big)
\nonumber\\
&=&\Big(G_{a\bar b;c}(\gamma,y)\dot\gamma^c
+G_{a\bar b;\bar c}(\gamma,y)\dot{\bar\gamma}^c
-G_{a\bar bc}(\gamma,y)\Gamma^c_{;e}(\gamma,y)\dot\gamma^e
+G_{a\bar b\bar c}(\gamma,y)\nonumber\\
&&\overline{\Gamma^c_{;e}(\gamma,y)\dot{\gamma}^e}\Big)
u^a\bar v^b
-G_{a\bar b}(\gamma,y)\Big(\Gamma^a_{;e}(\gamma,u)\dot\gamma^e
\bar v^b
+u^a\overline{\Gamma^b_{;e}(\gamma,v)\dot{\gamma}^e}\Big).
\label{affiso2}
\end{eqnarray}
Again by the complex Berwald condition,
\begin{equation}
\Gamma^a_{b;c}(\gamma,y)=\Gamma^a_{b;c}(\gamma)=
\Gamma^a_{b;c}(\gamma,u).\label{affiso3}
 \end{equation}
 By the formula $\Gamma^e_{a;c}=G^{\bar be}\delta_cG_{a\bar b}$,
expansion \eqref{affiso3}:
\begin{equation*}
\Big[G^{\bar be}\Big(G_{a\bar b;c}-\Gamma^d_{;c}G_{a\bar bd}\Big)\Big](\gamma,y)=
\Big[G^{\bar be}\Big(G_{a\bar b;c}-\Gamma^d_{;c}G_{a\bar bd}\Big)\Big](\gamma,u),
 \end{equation*}
thus
\begin{equation}
\Big(G_{a\bar b;c}-\Gamma^d_{;c}G_{a\bar bd}\Big)(\gamma,y)=
G_{e\bar b}(\gamma,y)\Big[G^{\bar se}\Big(G_{a\bar s;c}-\Gamma^d_{;c}G_{a\bar sd}\Big)\Big](\gamma,u).\label{affiso5}
 \end{equation}
By (1,1)-homogeneous of $G$,
 \begin{equation}
 G_{a\bar sd}(\gamma,u)u^a=0,\,
G_{a\bar s;c}(\gamma,u)u^a=G_{\bar s;c}(\gamma,u).\label{affiso6}
\end{equation}
By \eqref{affiso5} and \eqref{affiso6},
\begin{equation}
G_{a\bar b;c}(\gamma,y)u^a=
G_{e\bar b}(\gamma,y)\Gamma^{e}_{;c}(\gamma,u)
+\Gamma^d_{;c}(\gamma,y)G_{a\bar bd}(\gamma,y)u^a.\label{affiso7}
 \end{equation}
Take conjugation of \eqref{affiso7} and replace $u$ by $v$, we have
\begin{equation}
G_{\bar a b;\bar c}(\gamma,y)\bar v^a=
G_{\bar eb}(\gamma,y)\overline{\Gamma^{e}_{;c}(\gamma,v)}
+\overline{\Gamma^d_{;c}(\gamma,y)}G_{\bar ab\bar d}(\gamma,y)\bar v^a.\label{affiso8}
 \end{equation}
Substitute \eqref{affiso7} and \eqref{affiso8} into \eqref{affiso2}, we obtain $\dot g(t)\equiv0,$ thus $g(t)\equiv g(0).$

\textbf{Step 2.}
Notice $\gamma(0)=q,\,\gamma(1)=p,\,y(0),\,u(0),\,v(0)\in T^{1,0}_qM_1,\,u(1),\,v(1)\in T^{1,0}_1M_1$.
And $f$ is isometric at the point $p$, thus
\[g(0)=g_{y(0)}(u(0),v(0))=g(1)=g_{y(1)}(u(1), v(1))=g_{d f(y(1))}(df(u(1)), df(v(1))),\]
where
\[g_{d f(y)}(df(u), df(v))
=(F_2^2)_{a\bar b}(f(\gamma),df(y))df^a(u)\overline{df^\beta(v)}.\]
Now $y(t)$, $u(t)$, and $v(t)$ are results of parallel displacement of $y(0) ,\, u(0)$ and $v(0)$  along $\gamma,$ respectively.
Since $f$ is an affine transformation, transforming vectors that are parallel along $\gamma$ into vectors that are parallel along $f(\gamma)$, the vector fields $df(y(t))$, $df(u(t))$ and $df(v(t))$ are parallel displacement of $df(y(0))$, $df(u(0))$ and $df(v(0))$ along $f(\gamma)$, respectively.
 Thus, by Step 1, we have
\[g_{df(y(t))}(df(u(t)),df(v(t)))=g_{df(y(0))}(df(u(0)), df(v(0))).\]
Therefore
\[F_1^2(q,y(0))=g_{y(0)}(y(0), y(0))=g_{df(y(0))}(df(y(0)), df(y(0)))=F_2^2(f(q),df(y(0))).\]
Thus $f$ is an isometry.
\end{proof}

Now we give the characterization of strongly pseudoconvex K\"ahler--Berwald manifolds.
\begin{theorem}\label{characKB}
A connected, strongly pseudoconvex K\"ahler--Berwald manifold $(M,F)$ must be one of the  following four types:
\begin{enumerate}
\item [1)] K\"ahler space;

\item [2)] locally complex Minkowski space;

\item [3)] locally irreducible, locally symmetric, non-Hermitian, K\"ahler--Berwald manifold of rank $r \geq2$;

\item [4)] locally reducible, and in this case $(M, F)$ can be locally decomposed into a Descartes product of K\"ahler spaces, locally complex Minkowski spaces and locally irreducible, locally symmetric non-Hermitian K\"ahler--Berwald manifolds of rank $r \geq2.$
    \end{enumerate}
\end{theorem}
\begin{proof}Since Hermitian metrics and locally complex Minkowski metrics are K\"ahler--Berwald, we just need to consider when
$F$ is neither Hermitian nor locally complex Minkowski metric.
When $(M,F)$ is locally irreducible, we need to prove the third case holds.

\textbf{Step 1.} By Theorem \ref{Ai-thm}, the K\"ahler--Berwald metric $F$ induce a K\"ahler metric $h$, such that the nonlinear connection of $F$ coincides with the Chern connection of $h$, thus they have same unitary component holonomy group $Hol_p^o$.
By Theorem \ref{Berger}, when the manifold $M$ is simply connected and irreducible, then either $Hol_p^o$ acts transitively on the unit sphere with respect to the Hermitian metric $h$ in $T_pM$, or $(M,h)$ is a symmetric space of rank $\geq2$.
 The former case implies $F$ being a K\"ahler metric, since $F$ is also invariant under the action of  $Hol_p^o$ thus share same unit sphere with $h$, which contradict to our assumption. Thus $(M,h)$ is a Hermitian symmetric space of rank $\geq2$.

\textbf{Step 2.} For any point $p\in  M$, let  $\sigma_p$  be a holomorphically isometric geodesic symmetry of $(M,g)$ at $p$, with isolated fixed point $p$. By Lemma \ref{Chern--Finsler nonlinear-geod}, $F$ and $g$ have same geodesics, thus  $\sigma_p$  is a geodesic symmetry of $( M,F)$ at $p$, with isolated fixed point $p$.

\textbf{Step 3.} Below we prove $\sigma_p$ is also an isometry with respect to $F$.
By Lemma \ref{aff-iso}, this could be done if we can prove   $\sigma_p$ is an affinely transformation with respect to the Chern--Finsler nonlinear connection of $F$, and isometric at one point.
The former condition is a result of $\sigma_p$ being a holomorphically isometric geodesic symmetry of the K\"ahler metric $g$, whose Chern connection coincides with the Chern--Finsler nonlinear connection of $F$.
By {Step 2}, $\sigma_p$  is a locally biholomorphic geodesic symmetry of $( M,F)$ at $p$, with isolated fixed point $p$. Thus by Lemma \ref{Lemma}, $(d\sigma_p)\vert_p=-Id,$ hence $F(p,v)=F(p,-v)=F(\sigma_p(p),(d\sigma_p)\vert_p(v)),\;\forall \,v\in T^{1,0}_pM.$ So $\sigma_p$  is an isometry of $(M,F)$.

By Steps 1-3, $\sigma_p$ is a holomorphically isometric geodesic symmetry of $(M,F)$ at $p$ with $p$ as isolated fixed point, thus $F$ is a locally symmetric complex Finsler metric.

\textbf{Step 4.} Assume now $(M,F,\nabla)$ be locally reducible, and $\nabla$ the nonlinear connecction of $F$. Thus $F$ and the corresponding K\"ahler metric $h$ share same holonomy group $Hol_h$. For small simply connected neighborhood $U\subset M$, let $(U_1,h_1,\nabla^{(1)})\times\cdots (U_k,h_k,\nabla^{(k)})$ be the Descartes product of irreducible invariant K\"ahler submanifolds with respect to the product holonomy group $Hol_h^o=Hol^{(1)o}\times\cdots\times Hol^{(k)o}$ of $(U,h,\nabla)$.
Notice that $h(z,v)=\sum_jh\vert_{U_j}(z^{(j)},v^{(j)}),$ thus $(U_j, h\vert_{U_j},\nabla^{(j)})$ are totally geodesic submanifolds of $(U,h)$.

\textbf{Step 5.} By Proposition 4.1 of \cite{Ai95}, a complex Berwald metric is invariant under parallel displacement with respect to its nonlinear connection.
On each $U_j,$ consider the Minkowski norm $F^{(j)}_{p}$ induced by the restriction of $F$ at a fixed point $p\in U$, this norm is invariant under the action of $Hol_{p}^{(j)o}$, which are parallel displacements with respect to the Chern connection of $(M,h)$, thus with respect to the nonlinear connection of $(M,F)$.
Thus we can extend the norm $F^{(j)}_p$ to every point of $U_j$ through parallel translation with respect to $\nabla^{(j)}$, which coincides with the parallel translation with respect to $\nabla,$ since  $(U_j, h\vert_{U_j},\nabla^{(j)})$ are totally geodesic submanifold. The extension just gives the restriction of $F$ at each $U_j,$ and we obtain a strongly pseudoconvex complex Finsler metric $F^{(j)}$ on $U_j.$

We claim that $\nabla^{(j)}$ is the Chern--Finsler
nonlinear connection of $(U_j,F^{(j)})$.  Denote its connection coefficients by $\Gamma^{a_j}_{b_j;c_j},$ and consider the
corresponding connection on the vertical bundle
\[
\nabla^{(j)}: \Gamma(\mathcal V^{(j)}) \rightarrow \Gamma((T^{1,0} \tilde U^{(j)}))^\ast\times \mathcal V^{(j)})
\]
by
\[
V^{(j)}\mapsto (dV^{a_j}+V^{b_j}\Gamma^{a_j}_{c_j;d_j } dz^{c_j})\dot\partial_{a_j}.
\]
It is by definition a complex linear vertical connection.

For $u^{(j)},\, w^{(j)}\in \Gamma(\mathcal V^{(j)})$
and $X^{(j)}\in \mathcal H^{(j)},$ since $\pi^{(j)}_\ast\pi^{\ast0}(v^{(i)})=0,i\neq j$, we have
\[\langle u^{(j)}, w^{(j)} \rangle_{v^{(j)}}
=\langle \pi^{\ast0} u^{(j)}, \pi^{\ast0} w^{(j)} \rangle_{\pi^{\ast0} v^{(j)}}.\]
Thus
\begin{eqnarray}
X^{(j)}\langle u^{(j)}, w^{(j)} \rangle_{v^{(j)}}
=\pi^{\ast0} X^{(j)}\langle \pi^{\ast0} u^{(j)}, \pi^{\ast0} w^{(j)} \rangle_{\pi^{\ast0} v^{(j)}}.\label{dR-metric1}
 \end{eqnarray}
By the matrical property of the Chern--Finsler connection $\nabla$ of $(M,F)$, we have
\begin{eqnarray}
&&\pi^{\ast0} X^{(j)}\langle \pi^{\ast0} u^{(j)}, \pi^{\ast0} w^{(j)} \rangle_{\pi^{\ast0} v^{(j)}}\nonumber\\
&=&\langle \nabla_{\pi^{\ast0} X^{(j)}}\pi^{\ast0} u^{(j)}, \pi^{\ast0} w^{(j)} \rangle_{\pi^{\ast0} v^{(j)}}
+\langle \pi^{\ast0} u^{(j)}, \nabla_{\overline{\pi^{\ast0} X^{(j)}}} \pi^{\ast0} w^{(j)} \rangle_{\pi^{\ast0} v^{(j)}}\nonumber\\
&=&\langle \pi^{\ast0}(\nabla^{(j)}_{ X^{(j)}} u^{(j)}), \pi^{\ast0} w^{(j)} \rangle_{\pi^{\ast0} v^{(j)}}
+\langle \pi^{\ast0} u^{(j)}, \pi^{\ast0}(\nabla^{(j)}_{\overline{ X^{(j)}}} w^{(j)}) \rangle_{\pi^{\ast0} v^{(j)}}\nonumber\\
&=&\langle \nabla^{(j)}_{ X^{(j)}} u^{(j)}, w^{(j)} \rangle_{ v^{(j)}}
+\langle u^{(j)}, \nabla^{(j)}_{\overline{ X^{(j)}}} w^{(j)} \rangle_{ v^{(j)}}, \label{dR-metric2}
\end{eqnarray}
where in the second equality we viewed the fact that the Chern--Finsler nonlinear connection of $F$ is also the Chern connection of $h$, and $(U^{(j)},h\vert_{U^{(j)}})$ is a de Rahm factor complex submanifold.
By \eqref{dR-metric1} and \eqref{dR-metric2}, we have
\begin{eqnarray}
X^{(j)}\langle u^{(j)}, w^{(j)} \rangle_{v^{(j)}}
=\langle \nabla^{(j)}_{ X^{(j)}} u^{(j)}, w^{(j)} \rangle_{ v^{(j)}}
+\langle u^{(j)}, \nabla^{(j)}_{\overline{ X^{(j)}}} w^{(j)} \rangle_{ v^{(j)}}, \label{dR-metric3}
\end{eqnarray}
And \eqref{dR-metric3} implies $\Gamma^{a_j}_{b_j;c_j}$ is also the Chern--Finsler nonlinear connection coefficient of $F.$
Thus $(U^{(j)},F^{(j)})$ is a K\"ahler--Berwald manifold, and $(U^{(j)},F^{(j)},\nabla^{(j)}),j=1,\cdots,k,$ provide the Descartes decomposition of $(U,F,\nabla)$.

Since now $(U^{(j)},F^{(j)},\nabla^{(j)})$ is locally irreducible, thus $F^{(j)}$ is one of locally Minkowski metrics, K\"ahler metrics and non-Hermitian locally irreducible locally symmetric complex Finsler metrics.
\end{proof}

Now we prove the following de Rham docomposition theorem.
\begin{theorem}\label{KBDR}
A simply connected, complete, strongly pseudoconvex K\"ahler--Berwald manifold $(M,F)$ can be decomposed into the Descartes product of a complex Minkowski space,  simply connected complete irreducible K\"ahler spaces and  simply connected non-Hermitian complete irreducible symmetric K\"ahler--Berwald manifolds of rank $r \geq2.$
\end{theorem}
\begin{proof}
By Aikou\cite{Ai95}, the K\"ahler--Berwald metric $F$ induce a K\"ahler metric $h$, such that the Chern--Finsler nonlinear connection of $F$ coincides with the Chern connection of $h$, thus $(M,h)$ is a simply connected, complete, K\"ahler manifold.
By Theorem \ref{DEKahler}, $(M,h)$ is holomorphically isometric to the direct product
$(M^{(0)},h^{(0)})\times (M^{(1)},h^{(1)})\times\cdots\times (M^{(k)},h^{(k)})$,
where $(M^{(0)},h^{(0)})$ is a complex Euclidean space (possibly of dimension 0), and $(M^{(1)},h^{(1)}),\cdots,(M^{(k)},h^{(k)})$ are all simply connected, complete, irreducible K\"ahler manifolds.

Let $Hol^{(j)}$ be the holonomy group of $(M^{(j)},h^{(j)})$.
Notice $(M^{(j)},h^{(j)})$ are all totally geodesic complex submanifolds of $(M,h)$, hence also totally geodesic submanifolds of $(M,F)$. Both $h^{(j)}$ and $ F\vert_{M^{(j)}}$ are invariant under the action of  $Hol^{(j)o}$.
If  $Hol^{(j)o}$ acts transitively on the unit sphere of $T^{1,0}_{p^{(j)}}M^{(j)}$, the invariance of $F\vert_{M_k}$ under the action of $Hol^{(j)o}$ implies that $F\vert_{M^{(j)}}$ is a K\"ahler metric, hence $(M^{(j)},F\vert_{M^{(j)}})$ is a simply connected, complete, irreducible K\"ahler manifold.
If $Hol^{(j)o}$ acts not transitively, by Theorem \ref{Berger}, $(M^{(j)},h^{(j)})$ is a Hermitian symmetric manifold of rank $\geq2.$
Since $(M^{(j)},F\vert_{M^{(j)}})$ share same Chern--Finsler nonlinear connection with $h^{(j)}$, by same method as Steps 2-3 in the proof of Theorem \ref{characKB}, $(M^{(j)},F\vert_{M^{(j)}})$ is a symmetric complex Finsler space.
At last, for $(M^{(0)},h^{(0)})$, $h^{(0)}$ has trivial connection, so does $F\vert_{M^{(0)}},$ thus  $F\vert_{M^{(0)}}$  is a Minkowski metric.
\end{proof}

By Theorem \ref{thm1} and Theorem \ref{KBDR}, we have the following corollary.
\begin{theorem}\label{rank1}
There is no non-Hermitian, strongly pseudoconvex symmetric complex Finsler metric on a symmetric space of rank 1.
\end{theorem}

\begin{remark}
This result replaced the holomorphic invariance condition, required in the corresponding results of \cite{Zh2} and \cite{CGZ}, by assuming the symmetric complex Finsler metric condition.
\end{remark}

Using the de Rham decomposition theorem for K\"ahler--Berwald spaces, many classical results in K\"ahler geometry can be generalized to the K\"ahler--Berwald setting. For example, the generalized Frankel conjecture resolved by Mok \cite{Mok88} can be generalized as  follows.

\begin{theorem}\label{KBFrankel4}
Let $(M, F)$ be a compact, strongly pseudoconvex, K\"ahler--Berwald manifold of nonnegative holomorphic bisectional curvature.
Then  the universal cover $(\widetilde{M}, \widetilde{F})$ of $(M, F)$ can be decomposed into the Descartes product with  factors holomorphically isometric to:
\begin{enumerate}
  \item [1)]  complex Minkowski spaces;
  \item [2)]  complex projective spaces endowed with K\"ahler metrics of nonnegative holomorphic bisectional curvature;
  \item  [3)] irreducible compact-type symmetric complex Finsler spaces of rank $\geq2$.%, with each $F_l$ is unique on $M_l$ in the sense of Chern--Finsler nonlinear class.
\end{enumerate}
\end{theorem}
\begin{proof}
By Theorem \ref{KBDR}, $(\widetilde{M}, \widetilde{F})$ can be decomposed into the Descartes product of three kinds of manifolds.

 The first kind contains just a complex Minkowski space $(\mathbb C^{n_0},F_0)$.

  The second kind contains some irreducible compact simply connected complete K\"ahler manifolds $(M_j,g_j),\,j=1,\cdots,\widetilde k,$ and in this case, of nonnegative holomorphic bisectional curvature. Thus by the generalized Frankel conjecture proved by N. Mok \cite{Mok88}, $(M_j,g_j)$ is holomorphically isometric to either
  $(\mathbb{P}^{n_j},h_j),$ a complex projective space endowed with a K\"ahler metric of nonnegative holomorphic bisectional curvature, or $(\mathfrak G^{n_{l}},h_{l})$, an irreducible compact-type symmetric complex Hermitian spaces of rank $\geq2$, endowed with the canonical Hermitian symmetric metric;

 The last kind contains some irreducible compact-type symmetric complex Finsler spaces of rank $\geq2$ $(\mathfrak G^{n_{l}},F_{l}),\,l=k+1,\cdots,\widetilde s, $ and each $F_l$ is unique in the sense of sharing different Chern--Finsler nonlinear connection.
 \end{proof}

\section{Uniformization and rigidity theorems}\label{sec-unif}

In this section, we consider the uniformization theorem for strongly pseudoconvex K\"ahler--Berwald manifolds, and rigidity theorems for the Kobayashi and Carath\'eodory metrics.

\subsection{The uniformization theorem for K\"ahler--Berwald manifolds}\label{sub-cHSC}

\begin{theorem}\label{uniformization}
Every simply connected, complete K\"ahler--Berwald manifold of constant holomorphic sectional curvature $c$ is holomorphically isometric to

\begin{enumerate}
  \item  [1)] $\mathbb B^n$ endowed with the Bergman metric if  $c=-4$;
  \item  [2)] $\mathbb C^n$ endowed with a Minkowski metric  if $c=0$;
  \item  [3)] $\mathbb {P}^n$ endowed with the Fubini--Study metric if $c=4$.
\end{enumerate}\end{theorem}
\begin{proof}
\textbf{Step 1.}  By Theorem \ref{KBDR}, $M=M^{(0)}\times M^{(1)}\times\cdots\times M^{(k)},$ where $M^{(0)}=\mathbb C^{n_0}$, $M^{(i)},\,i=0,\cdots, k,$ are irreducible totally geodesic complex submanifolds of complex dimension $n_i$, endowed with the restriction metrics $F^{(i)}$ from $F$. The constant holomorphic sectional curvature condition gives
\begin{equation}
\begin{aligned}
cG^2&=-2G_\alpha\partial_{\bar \beta}\Gamma^\alpha_{\gamma;\sigma}v^\gamma v^\sigma\bar v^\beta\\
&=-2\sum_{i,j,k,l}
\sum_{\alpha_i}\sum_{\beta_j}\sum_{\gamma_k}\sum_{\sigma_l}
G_{\alpha_i}\partial_{\bar{\beta_j}}
\Gamma^{\alpha_i}_{\gamma_k;\sigma_l}
\bar v^{\beta_j}v^{\gamma_k} v^{\sigma_l}.
\end{aligned}
\label{total-cur''}
\end{equation}
Since the Chern--Finsler nonlinear connection of $F$ coincides with the Chern connection of a K\"ahler metric $h$,
then the de Rham decomposition theorem of $(M,h)$ gives
 \[\Gamma^{\alpha_i}_{\gamma_k;\sigma_l}
 =\delta^i_k\delta^i_l\Gamma^{\alpha_i}_{\gamma_k;\sigma_l},
 \]
which gives the connection coefficients of $h^{(i)}$ and $F^{(i)}$.
Thus \eqref{total-cur''} reduces to
\begin{equation}
cG^2
=-2\sum_{i}\sum_{\alpha_i}\sum_{\beta_i}\sum_{\gamma_i}\sum_{\sigma_i}
G_{\alpha_i}\partial_{\bar{\beta_i}}\Gamma^{\alpha_i}_{\gamma_i;\sigma_i}
\bar v^{\beta_i}v^{\gamma_i} v^{\sigma_i}.\label{total-cur'}
\end{equation}
Restricted at $\pi^{\ast0}(v^{(i)})=(0,\cdots,0,v^{(i)},0,\cdots,0)$, we have
\begin{equation}
cG^2(\pi^{\ast0}(v^{(i)}))
=-2\sum_{\alpha_i}\sum_{\beta_i}\sum_{\gamma_i}\sum_{\sigma_i}
G_{\alpha_i}(\pi^{\ast0}(v^{(i)}))
\partial_{\bar{\beta_i}}\Gamma^{\alpha_i}_{\gamma_i;\sigma_i}
\bar v^{\beta_i}v^{\gamma_i} v^{\sigma_i},
\end{equation}
that is, with each fixed $z^{(j)},\,j\neq i,$
\begin{equation}
c(G^{(i)}(v^{(i)}))^2
=-2\sum_{\alpha_i}\sum_{\beta_i}\sum_{\gamma_i}\sum_{\sigma_i}
G^{(i)}_{\alpha_i}(v^{(i)})
\partial_{\bar{\beta_i}}\Gamma^{\alpha_i}_{\gamma_i;\sigma_i}
\bar v^{\beta_i}v^{\gamma_i} v^{\sigma_i},\label{cons-curv-irre}
\end{equation}
which says that $F^{(i)}$ has constant holomorphic sectional curvature $c$.

\textbf{Step 2.} Since $(M^{(i)},F^{(i)})$ is a simply connected, complete, irreducible K\"ahler--Berwald manifold with constant holomorphic sectional curvature $c$, by Theorem \ref{characKB} and the uniformization theorem for K\"ahler manifolds, $(M^{(j)},F^{(j)})$ is one of the following:
\begin{itemize}
  \item [(i)] simply connected complete K\"ahler manifold with constant holomorphic sectional curvature $c$;
  \item [(ii)] $\mathbb C$ endowed with the Euclidean metric;
  \item [(iii)] non-Hermitian, irreducible, symmetric K\"ahler--Berwald manifold of rank $\geq2,$ and with constant holomorphic sectional curvature $c$.
\end{itemize}
In case of $(ii)$, $M=\mathbb C^n,$ and the Chern--Finsler nonlinear connection of $F$ vanishes. Thus
\[0=\Gamma^i_{k;j}v^k=\Gamma^i_{;j}
=G^{\bar li}G_{\bar l;j},
\]
thus
$G_{\bar l;j}=0$ since $(G^{\bar li})>0$, and $G_{;j}=G_{\bar l;j}\bar v^l=0.$ Thus the real valued map $G$ depends only on $v$, $F$ is a Minkowski metric.

\textbf{Step 3.} For case $(i)$,
%We claim that in case $(iii)$, the rank of $(M^{(i)},F^{(i)})$ must  be $1$.
By the uniformization theorem of K\"ahler manifolds, $(M^{(i)},F^{(i)})$ is either $\mathbb B^{n_i}$ endowed with the Bergman metric
or $\mathbb P^{n_i}$ endowed with the Fubini--Study metric, and $M$ is biholomorphic to either $\mathbb B^{n_1}\times \cdots \times \mathbb B^{n_k}$ or $\mathbb P^{n_1}\times \cdots \times \mathbb P^{n_k}$, depending on the signature of $c$. Thus, by the curvature of the Bergman metric of $\mathbb B^{n_i}$ or the Fubini--Study metric of $\mathbb P^{n_i}$,
\eqref{total-cur'} provides
\begin{equation}
\begin{aligned}
4G^2(0,v)
&=2\sum_{i}\sum_{\alpha_i}\sum_{\beta_i}\sum_{\gamma_i}\sum_{\sigma_i}
G_{\alpha_i}(0,v)
\left(\delta^{\alpha_i}_{\gamma_i}\delta_{\sigma_i\beta_i}
+\delta^{\alpha_i}_{\sigma_i}\delta_{\gamma_i\beta_i}\right)
\bar v^{\beta_i}v^{\gamma_i} v^{\sigma_i}\\
&=4\sum_{i}\sum_{\alpha_i}\sum_{\beta_i}\sum_{\sigma_i}
G_{\alpha_i}(0,v)v^{\alpha_i}
\delta_{\sigma_i\beta_i}\bar v^{\beta_i} v^{\sigma_i}\\
&=4\sum_{i}\sum_{\alpha_i}G_{\alpha_i}(0,v)v^{\alpha_i}\| v^{(i)}\|^2,
\end{aligned}\label{re-curv-1}
\end{equation}
where $0$ is the origin of $\mathbb B^{n_1}\times \cdots \times \mathbb B^{n_k}$ or the center of the affine chart $[1:z^{(1)1}:\cdots:z^{(1)n_1}]\times \cdots \times[1:z^{(k)1}:\cdots:z^{(k)n_k}]$
for $\mathbb P^{n_1}\times \cdots \times \mathbb P^{n_k}$.

Along a real curve $\xi: (a,b)\subset \mathbb R \to T^{1,0}_0M,$ denoting $g:=G(0,\cdot)$,
\begin{equation}
(g\circ \xi)'(t)
=2\operatorname{Re}\sum_{i}\sum_{\alpha_i}g_{\alpha_i}(\xi(t))
(\xi^{(i)\alpha_i})'(t).\label{der-met-curv}
\end{equation}
To use \eqref{re-curv-1}, we consider the solution for
\begin{equation}
2(\xi^{(i)\alpha_i})'(t)
=\xi^{(i)\alpha_i}(t)\|\xi^{(i)}(t)\|^2.\label{eq-curve}
\end{equation}
Notice the direction of $\xi'$ coincides with $\xi$, we can write $\xi^{(i)\alpha_i}(t)=r^{(i)\alpha_i}(t)
e^{\sqrt{-1}\theta^{(i)\alpha_i}},$ and \eqref{eq-curve} reduces to
$2{r^{(i)}}'(t)=r^{(i)}(t)\|r^{(i)}\|^2.$
Take $s:=\|r^{(i)}\|,$ then
\[
s'=\frac{\langle r^{(i)},r^{(i)'} \rangle}{\|r^{(i)}\|}
=\frac12 s^3.
\]
Thus $d(\frac{1}{s^2})=-1,$
$\frac{1}{s^2}=c-t.$ Evaluate $r^{(i)}(0)=v^{(i)}$ with fixed $v\neq0$, providing
\[
\|r^{(i)}\|^2=s^2=\frac{\vert v^{(i)}\vert^2}{1-t\vert v^{(i)}\vert^2}.
\]
This provides a solution of \eqref{eq-curve} by
\[\xi^{(i)}(t):=\frac{v^{(i)}}{\sqrt{1-t\|v^{(i)}\|^2}},\quad t<\frac1R,\,i=1,\cdots,k,\]
where $R:=\max\limits_i\|v^{(i)}\|^2.$
Thus \eqref{der-met-curv}, \eqref{eq-curve} and \eqref{re-curv-1} imply
\begin{equation}
(g\circ \xi)'(t)=g^2(\xi(t)).
\end{equation}
Thus
\begin{equation}
\left(\frac1{g\circ \xi}\right)'(t)=-1,
\end{equation}
and
\[
\int_0^t\left( \frac{1}{g\circ \xi}\right)'(t)dt
=\frac1{g(\xi(t))}-\frac1{g(v)}=-t.
\]
Notice
\begin{equation}
\frac1{g(v)}=\frac1{g(\xi(t))}+t\to\frac1R\quad (t\to \frac1R-),\label{g-lim}
\end{equation}
because of
 $\|\xi^{(i)}\|\to +\infty$ when
$t\to\frac1R-$, and $g(\xi)\ge \min\limits_{w\in\mathbb S^{2n}}g(w)\|\xi\|^2$ by homogeneity.
Thus the independence of  $g(v)$ with respect to $t$, and the limit in \eqref{g-lim} imply
\[g(v)=R=\max\limits_i\|v^{(i)}\|^2.\]
By the arbitrary choice of $v\neq0,$
\[g(v)=\max\limits_i\|v^{(i)}\|^2,\quad\forall\, 0\neq v\in T^{1,0}_0M.\]
Thus $g$ is smooth on $T^{1,0}_0M\setminus\{0\}$ if and only if
$k=1$, thus
\[g(v)=\|v\|^2,\quad\forall\, 0\neq v\in T^{1,0}_0M.\]
Since $G$ is symmetric, thus homogeneous, the expression of $g(\cdot)=G(0;\cdot)$ implies that $G$ is a K\"ahler metric. Thus the uniformization theorem for K\"ahler manifolds provides the conclusion.

\textbf{Step 4.} We prove that the case $(iii)$ is impossible by contradiction.
Suppose $(M^{(i)},F^{(i)})$ is an irreducible symmetric complex Finsler space of rank $r\geq 2$.
Let $h$ be the Hermitian symmetric metric whose Chern connection coincides with the Chern--Finsler nonlinear connection of $F^{(i)}$.
By the polydisc and polysphere theorems (cf. Theorem 1 in Chapter 5 of \cite{Mok}), when $M^{(i)}$ is of noncompact type,
there exists a totally geodesic complex submanifold $S\subset M^{(i)}$ such that $(S,h\vert_S)$ is isometric to the Poincar\'e polydisk $ (\mathbb D^r,\rho_o)$; and when $M$ is of noncompact type,
there exists a totally geodesic complex submanifold containing $S$ as an open subset such that $(S,h\vert_S)$ is isometric to the polysphere $((\mathbb P^1)^r,\rho_c)$ equipped with a product Fubini--Study metric.

Denote the restriction metric of $G^{(i)}$ on $S$ by $\hat G^{(i)}:=G^{(i)}\vert_{S}$,  by the Step 5 in proof of Theorem \ref{characKB}, $\hat G^{(i)}$ is a non-Hermitian symmetric K\"ahler--Berwald metric.
Since $S$ is totally geodesic complex submanifold of $M^{(i)}$,
the holomorphic sectional curvature of $\hat G^{(i)}$ coincides with that of $G^{(i)}$, which is thus also constant $c$.
Thus $S$ is either $\mathbb B^{n_i}$ or $\mathbb P^{n_i}$ by Steps $1-3$, since $S$ is a production of symmetric space of rank-$1$. Which is a contradiction.
\end{proof}

From the proof of Theorem \ref{characKB}, we see that each locally irreducible de Rham factor $(U^{(j)},F^{(j)})$ of a K\"ahler--Berwald manifold is a totally geodesic K\"ahler--Berwald sub-manifold.
And from Step 1 in the proof of Theorem \ref{uniformization}, we also obtain the coincidence for the $h\bar h$-part curvature between  $F^{(j)}$ and $F$. Thus we have the following results.
\begin{theorem}\label{curv-sub-KB}
For a K\"ahler--Berwald manifold $(M,F)$, let $(U^{(j)},F^{(j)}),j=1,\cdots,k,$ be its locally irreducible de Rham factors. Then  each $(U^{(j)},F^{(j)})$ is a totally geodesic K\"ahler--Berwald sub-manifold of $(M,F)$, and for fixed  $z^{(i)},\,i\neq j,$
\begin{align}
(R^{(j)})^{\alpha_j}_{\beta_j;\gamma_j,\overline{ \sigma j}}\left(z^{(j)},v^{(j)}\right)
&=R^{\alpha_j}_{\beta_j;\gamma_j,\overline{ \sigma j}}\left(z,\pi^{\ast0}(v^{(j)})\right),\\
(R^{(j)})_{\beta_j\overline{\mu_j};\gamma_j,\overline{ \sigma j}}\left(z^{(j)},v^{(j)}\right)
&=R_{\beta_j\overline{\mu_j};\gamma_j,\overline{ \sigma j}}\left(z,\pi^{\ast0}(v^{(j)})\right).
\end{align}
Especially,
\begin{align}
K^{(j)}\left(z^{(j)},v^{(j)}\right)
&= K\left(z,\pi^{\ast0}(v^{(j)})\right),\quad \forall\, v^{(j)}\in \tilde M^{(j)}_{z^{(j)}};\\
B^{(j)}\left(z^{(j)},v^{(j)},w^{(j)}\right)
&=B\left(z,\pi^{\ast0}(v^{(j)}), \pi^{\ast0}(w^{(j)})\right),\quad \forall\, v^{(j)},w^{(j)}\in \tilde M^{(j)}_{z^{(j)}}.
\end{align}
where \[\pi^{\ast0}(v^{(j)})
=\left(0,\cdots,0,v^{(j)},0,\cdots,0\right),\,
\pi^{\ast0}(w^{(j)})
=\left(0,\cdots,0,w^{(j)},0,\cdots,0\right).\]
\end{theorem}
%if $F$ has positive/semi-positive/negative/semi-negative  holomorphic sectional curvature, so do the metrics $F^{(i)},\,i=1,\cdots,k$. Thus many results in K\"ahler geometry can be generalized to the K\"ahler--Berwald setting.%, such as the generalized Frankel conjecture

In 2017, H. Xia and C. Zhong proposed the question \cite{XZ}:
is there any (strongly pseudoconvex) weakly K\"ahler--Finsler metric which is not K\"ahler--Finsler?
By Theorem \ref{thm:KF=KB}, K\"ahler--Finsler coincides with K\"ahler--Berwald, thus their question reduces to:
is there any strongly pseudoconvex weakly K\"ahler--Finsler metric which is not K\"ahler--Berwald?
This is a direct corollary by Theorem \ref{uniformization} and the Kobayashi metrics of smooth bounded strongly convex domains.
\begin{corollary}\label{non-Berwald}
Let $D\subset \mathbb C^n$ be a bounded strongly convex domain with smooth boundary, which is not biholomorphic to $\mathbb B^n$. Then its Kobayashi metric, is neither complex Berwald nor weakly complex Einstein--Finsler in the sense of Aikou.
\end{corollary}
\begin{proof}
By L. Lempert's result in \cite{Lem81}, the Kobayashi metric $F_K$ and the Carath\'eodory metric $F_C$ of $D$ coincide with each other, and they have constant holomorphic sectional curvature $-4$. And $(D,F_K)$ admits a unique complex geodesic through every $(z,v)\in T^{1,0}D$.

$F_K$ is a strongly %(strictly writed in their paper=with positive definite real Hessian=strictly convex Kobayashi indicatrix)
convex metric by Patrizio \cite{Pat}.
Thus by Corollary 3.4 of M. Abate and G. Patrizio \cite{AP2}, $F_K$ is a weakly K\"ahler--Finsler metric.

Since $D$ is not biholomorphic to $\mathbb B^n$, by Theorem \ref{uniformization}, $F_K$ is not K\"ahler--Berwald.
However, by \cite{XZ}, a (weakly) complex Berwald metric is weakly K\"ahler--Finsler if and only if it is K\"ahler--Finsler,
thus $F_K$ is not a complex Berwald metric.
Thus by Theorem \ref{SZ18EF}, $F$ is not weakly complex Einstein--Finsler.
\end{proof}

\begin{remark}
In \cite{CGZhou}, N. Cui, J. Guo and L. Zhou provided some singular weakly K\"ahler--Finsler metrics which are not K\"ahler--Finsler.
Smooth weakly K\"ahler--Finsler metric which is not K\"ahler--Berwald is still unknown before the present paper.
\end{remark}

\subsection{Rigidity theorems for the Kobayashi and Carath\'eodory metrics}\label{sub-KB-Kob}

%As corollaries of the former sections, we can generalize many results of K\"ahler metrics into K\"ahler--Berwald cases.
%For example, the uniformization theorem was obtained by results of Aikou \cite{Ai95} and Berger's result as Theorem \ref{Berger}.

In this section, we consider rigidity results for the Kobayashi metric $F_K$ and the Carath\'eodory metric $F_C$ under the K\"ahler--Berwald condition.

Recall that the definition of $F_K$ and $F_C$ of a complex manifold $M$ are defined, respectively, as follows:
\[
 F_K(z,v)=\inf\left\{\vert t\vert :t\in\mathbb C,
 f\in\operatorname{Hol}(\mathbb D,M),
 f(0)=z,\ tdf_0\left(\frac{\partial}{\partial \zeta}\right)=v\right\}
\]
and
\[
 F_C(z,v)=\sup\bigl\{\vert df_z(v)\vert :
 f\in\operatorname{Hol}(M,\mathbb D),\ f(z)=0\bigr\}.
\]

%Since Wong proved that any $C^2$-smooth bounded homogeneous domain is biholomorphic to the unit ball \cite{Wong77}, and the Bergman metric of a smooth bounded pseudoconvex domain is complete \cite{Oh81},
%thus there is also the conjecture of Cheng-Yau \cite{HHL26}:
%The Bergman metric of a $C^2$-smooth bounded  pseudoconvex domain $D$ is a Einstein metric if and only if $D$ is biholomorphic to the unit ball of same dimention.
%Cheng's conjecture was resolved recently by X. Huang and M. Xiao in \cite{Huang-Xiao}.
%, and C. Hsiao, X. Huang and X. Li proved the Cheng-Yau conjecture in case that the domain has real analytic boundary or the boundary is smooth in addition the h-extendable condition \cite{HHL26}.

\begin{theorem}\label{Kob-KB}
 For a complex manifold $M$, we have
\begin{itemize}
  \item [1)] $F_K$ is K\"ahler--Berwald if and only if it is K\"ahler;
  \item [2)] $F_C$ is K\"ahler--Berwald if and only if it is K\"ahler, when $M$ is in additionally assumed simply connected.
\end{itemize}
\end{theorem}
\begin{proof}
The sufficient direction is trivial.
We prove the necessary direction.

$1)$ Let $F_K$ be a strongly pseudoconvex K\"ahler--Berwald metric, and let $\pi:\widetilde M\rightarrow M$ be the holomorphic universal cover of $M$. Then $\pi$ is a locally biholomorphic, and locally isometric map with respect to the Kobayashi metrics of $\widetilde M$ and $ M$ (cf. Theorem 7.3.1 in \cite{GKK11}). Thus the Kobayshi metric $\widetilde{F_K}$ of $\widetilde M$ is also a strongly pseudoconvex K\"ahler--Berwald metric.

By the product property for the Kobayashi metric \cite{Jar-Pfl93}, if $\widetilde M=\widetilde M_1\times\widetilde M_2$ is a product manifold, $\widetilde{F_K}$ is defined by taking the maximum of $\widetilde{F_{K,\widetilde M_1}}$ and $\widetilde{F_{K,\widetilde M_2}}$, which is not strongly pseudoconvex.
Hence by Theorem \ref{characKB}, if $\widetilde{F_K}$ is not a K\"ahler metric, $\widetilde M$ is an irreducible Hermitian symmetric space of rank $\geq2$.
The later is impossible: $\widetilde {F_K}$ vanishes if $\widetilde M$ is an irreducible Hermitian symmetric space of compact type; $\widetilde{F_K}$ is also irregular if $\widetilde M$ is an irreducible Hermitian symmetric space of noncompact type with rank $\geq2$ (in this case, $\widetilde M$ is biholomorphic to an irreducible bounded symmetric domain, and $\widetilde{F_K}$ is defined also by taking maximum for corresponding eigenvalues).
Thus $\widetilde{F_K}$ is a K\"ahler metric, so does $F_K$.

$2)$ Assume $M$ is simply connected and $F_C$ is K\"ahler--Berwald.
Again by the product property for the Carath\'eodory metric \cite{Jar-Pfl89,Jar-Pfl99}, $M$ is irreducible. By  Theorem \ref{characKB}, if $F_C$ is not a K\"ahler metric, $M$ is an irreducible Hermitian symmetric space of rank $\geq2$.
The later case will not appear:
since $F_C$ vanishes on Hermitian symmetric space of compact type; and $F_C=F_K$ on Hermitian symmetric space of noncompact type, which is irregular if $M$ has rank $\geq2$.
\end{proof}

Below we give the proof for Corollaries \ref{1Kob-Ball} and \ref{1DWW-KB},
the rigidity for the Kobayashi or Carath\'eodory metric has locally K\"ahler--Berwald condition.
\begin{corollary}\label{KB-Ball}
Let $\Omega\subset \mathbb C^n,n\geq2,$ be a bounded strongly pseudoconvex domain with $C^2$ boundary, and assume that its Kobayashi metric $F_K$ is strongly pseudoconvex.
Then the following are equivalent:
\begin{enumerate}
    \item [1)] the universal cover of $\Omega$ is biholomorphic to $\mathbb B^n$;
  \item  [2)] $F_K$ is K\"ahler--Einstein and of constant holomorphic sectional curvature;
  \item  [3)] $F_K$ is weakly K\"ahler--Finsler, weakly complex Einstein--Finsler on $U \cap \Omega$, where $U\supseteq \partial \Omega$ is an open set;
  \item  [4)] $F_K$ is K\"ahler--Berwald on $U \cap \Omega$, where $U\supseteq \partial \Omega$ is an open set.
\end{enumerate}
\end{corollary}
\begin{proof}
Notice that the Kobayashi metric of $\mathbb B^n$ is K\"ahler--Einstein metric of constant holomotphic sectional curvature, and holomorphic covering map is isometric with respect to the Kobayashi metrics. Thus $(1)\Leftrightarrow (2)$ is due to Gaussier and Zimmer's result represented as Theorem \ref{Gaussier-Zimmer}.

 $(2)\Rightarrow (3)$ is by definition.
$(3)\Rightarrow(4)$ follows from Corollary \ref{SXZ},
and $(4)\Rightarrow (2)$ is by Theorem \ref{Kob-KB} and Theorem \ref{Gaussier-Zimmer}.
\end{proof}

%\begin{proof} Let $\pi:\widetilde D\rightarrow D$ be the universal cover of $D$, by condition  $\widetilde D$ is not a product doamin, thus it is irreducible. Since $\pi$ is isometry with respect to the Kobayashi metrics of $\widetilde D$ and $D$, the Kobayashi metrics $\widetilde{F_K}$ of $\widetilde D$ is also a strongly pseudoconvex K\"ahler--Berwald metric.
%If $\widetilde D$ is symmetric, then $\widetilde{F_K}$ is strongly pseudoconvex only if $\widetilde D$ is the unit ball, and we have done. Suppose $\widetilde D$ is not symmetric, thus by Theorem \ref{KBDR}, $\widetilde{F_K}$ is a K\"ahler metric, thus the results follows from Theorem \ref{Gaussier-Zimmer}.
%\end{proof}

%By Corollary \ref{non-Berwald} and Theorem \ref{KB-Ball}, we have
% \begin{theorem}\label{KFE-Ball2}
%Let $D\subset \mathbb C^n$ be a smooth bounded strongly convex domain, and $F_K$ its Kobayashi metric. Then $F_K$ is weakly complex Einstein--Finsler if and only if $D$ is biholomorphic to the unit ball.
%\end{theorem}

%In \cite{Gau-Zim}, H. Gaussier and A. Zimmer's proof of Theorem \ref{Gaussier-Zimmer} used an extendable property of Hermitian locally symmetric metrics, which is a metric analogue of the Hartogs theorem for holomorphic functions. Such extension theorem also holds for locally symmetric complex Finsler metrics.
Corollary \ref{1DWW-KB} can be proved similarly.
Besides, we can also using the following extendible property for locally symmetric complex Finsler metrics to provide another proof. Such extendible property is based on Gaussier--Zimmer's extendible property for locally symmetric Hermitian metrics \cite{Gau-Zim}.
\begin{lemma}\label{KB-ext}
  Suppose that $X$ is a Stein manifold with $dim_C \,X\geq2$, $K \subset\subset X$ is a compact subset, $X\setminus K$ is connected, and $G_{0}$ is a strongly pseudoconvex complex Finsler  metric on $X\setminus K$ which is complete at infinity. If $G_{0}$ is a locally symmetric complex Finsler metric, then there exists a complete locally symmetric complex Finsler metric $G$ on $X$ such that $G=G_{0}$ on $X \setminus K.$
\end{lemma}
\begin{proof}
$G_0$ is locally symmetric complex Finsler metric, by Theorem \ref{thm1}, $G_0$ is K\"ahler--Berwald metric, and its Chern--Finsler nonlinear connection coincides with the Chern connection of a locally symmetric Hermitian metric $g_0$ on $X\setminus K$.

By definition, for any smooth curve $c:[0,1]\rightarrow X\setminus K$, with $c(0)=p,\, c(1)=z,$ $g_0$ satisfies
\[ g_0(z,v)=g_0(p,\pi_c^{-1}(v)),\]
where $\pi_c$ is the parallel translation along $c$ with respect to the Chern connection of $g_0$.
By Theorem 1.1 of Gaussier and Zimmer in \cite{Gau-Zim}, there is a
complete locally symmetric Hermitian metric $g$ on $X$, such that
$g\vert_{X\setminus K}=g_0.$

For any $q\in K,$ there are locally holomorphically isometric geodesic symmetries $\sigma_i,i=1,\,\cdots,N,$ such that $p=\sigma_N\circ\sigma_{N-1}\circ \cdots\sigma_1(q).$ Without loss of generality we can assume $p=\sigma(q),$ $\sigma $ is a holomorphic isometric locally symmetry at $p_0\in X\setminus K.$
Hence $g(q;v)=g(\sigma(q),d\sigma(v))=g(p;d\sigma(v)).$
Define $G(q;v):=G_0(p;d\sigma(v))$, it's easy to see  $G(q;v)$ is independent of the choice of $\sigma$ since $g_0$ and $G_0$ share same geodesics near $p$. Thus we obtain a complex Finsler metric $G$ on $X$. By the proof of Theorem \ref{characKB}, holomorphically isometric geodesic symmetries of $g$ are also holomorphically isometric geodesic symmetries of $G$, thus $G$ is a complete locally symmetric complex Finsler metric.
\end{proof}

By Theorem \ref{uniformization},  Lemma \ref{KB-ext}  and Theorem \ref{DWW}, we have the following corollary.
%Since the universal cover mapping is not an isometry with resppect to the Carath\'eodory metrics, there may be no similar rigidity for the Carath\'eodory metric to be K\"ahler--Berwald. And Theorem \ref{DWW} can be generalized to K\"ahler--Berwald circumstance in addition simply connected condition.
\begin{corollary}\label{DWW-KB}
Let $\Omega \subset \mathbb{C}^{n},$ $n \geq2$, be a simply connected, bounded strongly pseudoconvex domain with $C^{2}$-boundary. Then the following are equivalent:
\begin{enumerate}
  \item $\Omega$ is biholomorphic to $\mathbb B^n$;
  \item The Carath\'eodory metric of $\Omega$ is a strongly pseudoconvex K\"ahler--Berwald metric on $U \cap \Omega$, where $U\supseteq \partial \Omega$ is an open set;
  \item The Carath\'eodory metric of $\Omega$ is a strongly pseudoconvex weakly K\"ahler--Finsler, weakly complex Einstein--Finsler metric  on $U \cap \Omega$, where $U\supseteq \partial \Omega$ is an open set.
\end{enumerate}
\end{corollary}
\begin{proof}
The directions $(1)\Rightarrow(2)$ and $(1)\Rightarrow(3)$ are evident,
and the direction $(3)\Rightarrow(2)$ is obtained by Corollary \ref{SXZ}.
We just need to prove $(2) \Rightarrow (1)$.
%Using Theorem \ref{Gau-Zim-nonempty}, in the proof of Theorem \ref{DWW}, Dong, Wang and Wong proved that: for all $z\in \Omega$ such that  $\delta_{\Omega}(z)<\epsilon,$ then there exists a nonempty open set $E_{\Omega}(z)$ of nonzero tangent vectors at $z$ where the Kabayashi metric and the Carath\'eodory metric agree, and the holomorphic sectional curvature of the Carath\'eodory metric $F_C$ is equad to $-2$ on $D\setminus \Sigma$, where $\Sigma:=\{z\in \Omega: \delta_{\Omega}(z)\geq \epsilon\}.$

% Denote by $\pi:\hat\Omega\rightarrow (\Omega\setminus \Sigma,F_C)$ the universal cover of $D\setminus \Sigma$. Notice that $F_C$ is K\"ahler--Berwald by condition, the pull back metric $\hat{F}_C$ of $F_C$ is K\"ahler--Berwald and of constant holomorphic sectional curvature $-2,$ thus by Corollary \ref{uniformization}, $\hat F_C$ is a constant multiple of the Bergman metric of the unit ball up to a biholomorphic isometry.
% Thus $F_C$ is a Hermitian locally symmetric metric restricted on $D\setminus \Sigma,$ thus the results follows from Theorem \ref{DWW}.

%Let $\pi:\widetilde D\rightarrow D$ be the universal cover of $D$,
Suppose (2) holds.
Since the Carath\'eodory metric $F_C$ is a strongly pseudoconvex metric, by the product property of $F_C$ \cite{Jar-Pfl89,Jar-Pfl99}, $\Omega$ is not a product domain. Thus $\Omega$ is locally irreducible since it is simply connected.
Notice that, by Theorem \ref{KBDR}, $F_C$ is either a K\"ahler metric or it is a non-Hermitian locally symmetric complex Finsler metric on $U \cap \Omega$. If $F_C$ is K\"ahler on $U \cap \Omega$, the results follows from Theorem \ref{DWW}.

Assume that $F_C$ is a non-Hermitian locally symmetric complex Finsler metric on $U \cap \Omega$. By the work of Graham (Proposition 5 of \cite{Graham}), the Carath\'eodory metric is complete on $\Omega$.
By Lemma \ref{KB-ext}, there exists a locally symmetric complex Finsler metric $\hat F_C$ on $\Omega$ such that
\[\hat F_C\vert_{U \cap \Omega}=F_C\vert_{U \cap \Omega}.\]

Notice that $\Omega$ is simply connected, $\hat F_C$  is a symmetric complex Finsler metric and $M$ is an irreducible bounded symmetric domain. Since $F_C$ is strongly pseudoconvex K\"ahler--Berwald metric on $U \cap \Omega$, and the only bounded symmetric domain
whose Carath\'eodory metric is strongly pseudoconvex is $\mathbb B^n$.
%And by the proof of Theorem \ref{DWW} in \cite{Dong-Wang-Wong}, $F_C$ has constant holomorphic sectional curvature $-2$ on $U\cap\Omega,$ so does $\widetilde F_C$. Since $\widetilde F_C$ is a symmetric metric, it is a K\"ahler--Berwald metric with constant holomorphic sectional curvature $-2$ on $D.$  Thus by Theorem \ref{uniformization}, $D$ is biholomorphic to the unit ball $\mathbb B^n.$
%$D$ is an irreducible bounded symmetric domain  of rank $\geq2$. Notice that the Carath\'eodory metric on such a domain is not differential, which contradict to the condition that $F_C$ is a strongly pseudoconvex metric. Thus $F_C$ is a K\"ahler metric.
%Since $\pi$ is isometry with respect to the Kobayashi metrics of $\widetilde D$ and $D$, the Kobayashi metrics $\widetilde{F_K}$ of $\widetilde D$ is also a strongly pseudoconvex K\"ahler--Berwald metric. If $\widetilde D$ is symmetric, then $\widetilde{F_K}$ is strongly pseudoconvex only if $\widetilde D$ is the unit ball, and we have done. Suppose $\widetilde D$ is not symmetric, thus by Theorem \ref{KBDR}, $\widetilde{F_K}$ is a K\"ahler metric, thus the results follows from Theorem \ref{Gaussier-Zimmer}.
\end{proof}

%A nature question is as follows.
%\begin{question}
%Does the simply connected condition of Theorem \ref{DWW-KB} can be removed?
%\end{question}

Another corollary of Theorem \ref{Kob-KB} is as follows.
\begin{corollary}\label{rig-unb-con-Kob}
For every convex domain $\Omega\subset\mathbb C^n$, its Kobayashi metric is K\"ahler--Berwald if and only if $\Omega$ is biholomorphic to $\mathbb B^n$.
\end{corollary}
\begin{proof}
The Kobayashi metric $F_K$ of a bounded convex domain coincides with its Carat\'eodory metric by Lempert \cite{Lem81}, and the unbounded convex domain case by \cite[Lemma 3.1]{Bracci-Saracco}. Thus the holomorphic sectional curvature of $F_K$ is $-4$.

If  $F_K$ is K\"ahler--Berwald, by Theorem \ref{Kob-KB}, it is K\"ahler, thus $\Omega$ is biholomorphic to $\mathbb B^n$ by the uniformization theorem for K\"ahler manifolds.
\end{proof}
\begin{remark}
Since every bounded homogeneous domain $\Omega\subset \mathbb C^n$ is biholomorphic to a affinely homogeneous Siegel domain by the Vinberg--Gindikin--Piatetski--Shapiro theorem \cite{VGPS63}, thus its Kobayashi metric is K\"ahler--Berwald if and only if $\Omega$ is biholomorphic to $\mathbb B^n$.
\end{remark}

\begin{remark}
The first-named author announced most of the core results in the K\"ahler--Berwald setting at the Young Mathematicians Workshop on Several Complex Variables, held in Incheon, Korea, in August 2026. The full proofs, as well as the equivalence between the K\"ahler--Berwald and K\"ahler--Finsler conditions, were finalized subsequently.

After the written of the present paper, we find that
J. Li, L. Li and C. Qiu also obtained the following results by different method \cite{LLQ26}: the uniformization theorem for K\"ahler--Finsler manifolds under a curvature symmetry condition \eqref{curv-sy},
and the K\"ahler--Finsler Kobayashi rigidity theorem for smooth bounded strongly pseudoconvex domains.
\end{remark}

\vskip0.4cm
\noindent
\textbf{Acknowledgement.}\ {\small
The second author was supported by the National Key R\&D Program of China (No. 2021YFA1003100) and the National Natural Science Foundation of China (Grant No. 12288201).

%\vskip0.4cm
%\noindent
%\textbf{Declarations}\\
%\noindent
%\textbf{Funding} (This work is partially supported by National Natural Science Foundation of China (Grant Nos. 12071386, 11671330)\\
%\noindent
%\textbf{Conflicts of interest/Competing interests} (There is no conflicts of interest)\\
%\noindent
%\textbf{Availability of data and material} (My manuscript has no associated data)\\
%\noindent
%\textbf{Code availability} (Not applicable)\\


%\input sn-sample-bib.tex%

\end{document}